\documentclass[11pt]{article}
\usepackage[utf8]{inputenc}
\usepackage[english]{babel}
\usepackage[hidelinks]{hyperref}
\usepackage{graphicx}
\usepackage{amsmath}
\usepackage{amssymb}
\usepackage{amsthm}
\usepackage{tikz-cd}
\usepackage{geometry}
\usepackage{xcolor}
\usepackage{etoolbox}
\usepackage[shortlabels]{enumitem}
\usepackage[center,md]{titlesec}
\patchcmd{\thebibliography}{\section*}{\section}{}{}
\theoremstyle{definition}
\newtheorem{thm}{Theorem}[subsection]
\newtheorem{cor}[thm]{Corollary}
\newtheorem{lem}[thm]{Lemma}
\newtheorem{prop}[thm]{Proposition}
\newtheorem{deff}[thm]{Definition}

\newtheorem{obs}[thm]{Remark}
\newcommand{\N}[0]{\mathbb{N}}
\newcommand{\Z}[0]{\mathbb{Z}}

\newcommand{\J}[0]{\mathcal{J}}

\begin{document}
\begin{center}
    \Large{\textbf{ON THE ABELIANIZATION OF CERTAIN GROUPS OF FORMAL POWER SERIES}}\vspace{0.8cm}\\
    \large{JAVIER PAVEZ CORNEJO\\
   email \href{mailto:javier.pavez.c@ug.uchile.cl}{javier.pavez.c@ug.uchile.cl}}\vspace{1cm}
\end{center}
\begin{abstract}
    We compute the abelianization of the Jennings group $\J_k(\Z)$ of powers series with constant coefficient $0$, linear coefficient equal to $1$ and vanishing coefficients in orders greater or equal than $2$ and less than $k$, where $k\geqslant2$. This is accomplished by directly dealing with the equivalence classes in the corresponding abelianizations, in contrast with the work of I. K. Babenko and S. A. Bogaty\u{\i} \cite{MR2413648}, who give an explicit abelianization morphism for the case $k=2$.
\end{abstract}
\tableofcontents

\newpage

\section{Introduction}

Let $R$ be a commutative ring. We define the Jennings group over $R$ to be 
\[\J(R):=\left\{f\in R[[x]]~\bigg|~f=x+\sum_{i\geqslant2}\alpha_ix^i\right\}.\]
In this generality, this group was first considered by S.A. Jennings (hence the name) in \cite{MR0061610}.

For each $k\geqslant2$, we consider the subgroup 
\[\J_k(R):=\left\{f\in\J(R)~\bigg|~f=x+\sum_{i\geqslant k}\alpha_ix^i\right\}.\]
The family of these subgroups will be the central object in this article, and we will get more acquainted with them later.
\\

The group $\J(R)$ has a natural topological group structure. This can be seen by identifying $R[[x]]$ with $R^\N$ as sets, and thus endowing it with the product topology (if $R$ is not explicitly said to be a topological ring, it can be assumed to be one with the discrete topology). Since the (composition) operation in $\J(R)$ can be encoded by sums and products of coefficients, this makes it a continuous operation in the subspace topology.
\\

As the naming convention would suggest, we are not the first studying such groups. Even though the study of groups of formal power series dates back to the works of S. Bochner and W.T. Martin \cite{MR0027863}, and M. G\^{o}to \cite{MR38351}, it is in his work \cite{MR0061610} that S.A. Jennings proves that $\J(R)$ is a topological group under composition, but not necessarily with the same topology that we defined. In his work, Jennings established many elementary results, one of which will be useful for us:

\prop\label{prop:cont}  One has $[\J_k(R),\J_k(R)]\subseteq\J_{2k}(R)$  for all $k \geq 2$.
\\\\
It is important to note that, in this context, the symbol $[\J_k(R),\J_k(R)]$ refers to the closure of the subgroup generated by the commutators of $\J_k(R)$, as opposed to simply such subgroup. Even though we stated that the topology studied by Jennings is not necessarily the same topology we will study in this article, this propositions still holds for our cases.

We will focus on the case of $R=\Z$. This presents some difficulties, in general by the fact that $\J(\Z)$ is homeomorphic to $\Z^\N$, which is not locally compact and restricts us from using results derived from the existence of a Haar measure, and in particular by the fact that $\Z$ has noninvertable elements.
\\

In \cite{MR2413648}, I. K. Babenko and S. A. Bogaty\u{\i}, among many other results, prove that the abelianization of $\J(\Z)$ is isomorphic to $\Z^2\bigoplus(\Z/2\Z)^2$. This was done by proving that the function

\begin{align*}
    \varphi:\mathcal{J}(\Z)&\to\Z\oplus\Z\oplus\Z/2\Z\oplus\Z/2\Z\\
    x+\sum_{i\geqslant2}\alpha_ix^i&\mapsto\left(\alpha_2,\alpha_3-\alpha_2^{~2},\alpha_4+\alpha_5+\frac{\alpha_3(\alpha_3+1)}{2}\pmod{2},\alpha_7+\alpha_5\alpha_3+\alpha_5\pmod{2}\right)
\end{align*}
is a surjective morphism whose kernel is $[\J(\Z),\J(\Z)]$. In contrast with this result, we will compute the abelianization of the groups $\J_k(\Z)$ for $k\geqslant2$ without an explicit abelianization morphism. This will be done by directly computing the equivalence classes in $H_1(\J_k(\Z))$. To make this possible, we will prove that for all $k\geqslant2$, there exists $c_k>k$ such that $\J_{c_k}(\Z)\subset [\J_k(\Z),\J_k(\Z)]$.  This will give us a simple method to verify some of the properties of the equivalence classes of $H_1(\J_k(\Z))$.
\\

Below we state our main result:
\begin{thm}\label{thm:princ}
    Let $k\geqslant2$.
\begin{enumerate}[$a)$]
\item  If $k\equiv 3\pmod{4}$, then \[H_1(\J_k(\Z))\cong\Z^k\bigoplus (\Z/2\Z)^\frac{k+1}{2}.\]

\item If $k\equiv2\pmod{4}$, then \[H_1(\J_k(\Z))\cong\Z^k\bigoplus (\Z/2\Z)^\frac{k+2}{2}.\]

\item If $k\equiv 1\pmod{4}$, then \[H_1(\J_k(\Z))\cong\Z^k\bigoplus\Z/4\Z\bigoplus (\Z/2\Z)^{\frac{k-3}{2}}.\]

\item If $k\equiv 0\pmod{4}$, then \[H_1(\J_k(\Z))\cong\Z^k\bigoplus\Z/4\Z\bigoplus (\Z/2\Z)^{\frac{k-2}{2}}.\]
\end{enumerate}
\end{thm}
~\\
As we mentioned before, this theorem does not give explicit abelianization morphisms. We hope to detect such morphisms in a future work.

\section{Preliminaries}

\subsection{Jennings Groups}

Before getting into the new results, we will familiarize ourselves with the Jennings groups. We will try to maintain the highest possible generality through this article. With this goal in mind, a torsion-free domain will be any integral domain $R$ such that the unique ring morphism $\Z\to R$ is injective. For this section, the ring $R$ will be a unitary commutative ring, unless specified otherwise.
\\

In general, finding the coefficients of $f\circ g$, for any $f,g\in\J(R)$, is not an easy task. Thankfully for us, given some conditions, we can give a formula for the first two non-linear coefficients of both $f\circ g$ and, in some cases, for $[f,g]$.

\prop \label{prop:formulas} Let $f=x+\alpha_mx^m+\alpha_{m+s}x^{m+s}+\cdots$ and
$g=x+\beta_nx^n+\beta_{n+r}x^{n+r}+\cdots$ be elements of $\J(R)$. Let $t=\min\{s,r,n-1\}$ and 
\[\gamma_1=\left\{\begin{array}{cc}
    (m+s)\alpha_{m+s}\beta_n & \text{if }t=s \\
   0  & \text{if }t\neq s
\end{array}\right.,~\gamma_2=\left\{\begin{array}{cc}
    m\alpha_m\beta_{n+r} & \text{if }t=r \\
   0  & \text{if }t\neq r
\end{array}\right.,~\gamma_3=\left\{\begin{array}{cc}
    \binom{m}{2}\alpha_m\beta_n^{~2} & \text{if }t=n-1 \\
   0  & \text{if }t\neq n-1
\end{array}\right. .\]

Then:
\begin{enumerate}[$a)$]
    \item There exists $E_{f,g}\in x^{m+n+t}R[[x]]$ such that
    \[f\circ g=f+g-x+m\alpha_m\beta_nx^{m+n-1}+(\gamma_1+\gamma_2+\gamma_3)x^{m+n+t-1}+E_{f,g}.\]
    
    \item If $m>n$ and  $r>s>n$, then \[[f,g]=x+(m-n)\alpha_m\beta_nx^{m+n-1}+C(m,n)\alpha_m\beta_n^2x^{m+2(n-1)}+\cdots\] where $C(m,n)=\binom{m}{2}-(m-n)(m+n-1)$.
\end{enumerate}

\begin{proof}~
    \begin{enumerate}[$1)$]
        \item See Proposition 2.3 of \cite{MR2413648}. The proof therein is done for integer coefficients, but it easily generalizes to $R$.

        \item Using the formula in $a)$ we have 
        \begin{align*}
            f\circ g-g\circ f&= (m-n)\alpha_m\beta_nx^{m+n-1}+\binom{m}{2}\alpha_m\beta_n^{~2}x^{m+2(n-1)}+\cdots\\
           \Leftrightarrow [f,g]-x&= (f\circ g-g\circ f)\circ(g\circ f)^{-1}\\
           &=(m-n)\alpha_m\beta_nx^{m+n-1}+C\alpha_m\beta_n^{~2}x^{m+2(n-1)}+\cdots\\
           \Leftrightarrow [f,g]&= x+(m-n)\alpha_m\beta_nx^{m+n-1}+C\alpha_m\beta_n^{~2}x^{m+2(n-1)}+\cdots.
        \end{align*}
        It follows from the previous equation that
        \begin{align*}
            [f,g]\circ(g\circ f)&=g\circ f+(m-n)\alpha_m\beta_n(g\circ f)^{m+n-1}+C\alpha_m\beta_n^{~2}(g\circ f)+\cdots\\
            &=g+f-x+n\beta_n\alpha_mx^{m+n-1}+(m-n)\alpha_m\beta_n(x+\beta_nx^n+\cdots)^{m+n-1}+C\alpha_m\beta_n^{~2}(x+\cdots)+\cdots\\
            &=g+f-x+((n+(m-n))\alpha_m\beta_n)x^{m+n-1}+((m+n-1)(m-n)+C)\alpha_m\beta_n^{~2}x^{m+2(n-1)}+\cdots.
        \end{align*}
       Since $[f,g]\circ g\circ f=f\circ g$, we can identify coefficients using the formula in $a)$. From this, we get the following:
        \[\binom{m}{2}=C+(m-n)(m+n-1).\]
        Thus $C$ depends only on $m$ and $n$, and equals $\binom{m}{2}-(m-n)(m+n-1)$.
    \end{enumerate}
    \vspace{-0.75cm}
\end{proof}

In some cases, we will simplify the formula for $f\circ g$ by writing

\[f\circ g=f+g-x+E,\]
with $E\in x^{m+n-1}R[[x]]$.\\

It is worth mentioning that the $l^{th}$ coefficient of $E_{f,g}$ is a multiple variable polynomial evaluated in $\{\alpha_i\}_{i=m}^{l-1}\cup\{\beta_i\}_{i=n}^{l-1}$.
\\

Another concept that will help us in this article is the following:

\deff Let $f=x+\sum_{i\geqslant2}\alpha_ix^i\in\J(R)\setminus\{x\}$. We denote the level of $f$ as
\[lvl(f):=\min\{i\geqslant2|\alpha_i\neq0\}.\]
For the special case of the identity element, we define $lvl(x)=\infty$.
\\

Said plainly, the level of an element of $\J(R)$ is the minimal power of $x$ that is larger than $1$ and whose coefficient is not $0$.
\\

Thanks to Proposition \ref{prop:formulas}, we deduce the following:
\cor \label{cor:nvl} Let $f$ and $g$ be elements in $\J(R)$. Then:
\begin{enumerate}[$a)$]
    \item $ lvl(f\circ g)\geqslant\min\{lvl(f),lvl(g)\}$. 
    
    \item $lvl(f^{-1})=lvl(f)$, and there exists $E_f\in x^{2lvl(f)-1}R[[x]]$ such that $f^{-1}=2x-f+R$.

    \item If $R$ is a torsion-free domain, then $lvl(f^l)=lvl(f)$. Furthermore, for all $l\in\Z^+$ there exists $E_l\in x^{2lvl(f)-1}R[[x]]$ such that $f^l=lf-(l-1)x+R_l$.

    \item $lvl([f,g])\geqslant lvl(f)+lvl(g)-1$. 
\end{enumerate}

This way, we can rewrite
\[\J_k(R)=\{f\in\J(R)|lvl(f)\geqslant k\}.\]
In particular $\J_2(\Z)=\J(\Z)$. We now take time to address the fact that, in the literature, these subgroups are usually indexed by $k-1$. We are breaking with this convention because it results in a more natural way of working with the subgroups and their elements, and some constants will become easier to manage. The maths underlying this work is nonetheless invariant, so for anyone who wishes to maintain convention, it suffices to translate where it should correspond.
\\

With this notation, it is easy to prove the following:

\prop Given $k\geqslant2$ and another integer $l$, one has $l\geqslant k$ if and only if $\J_l(R)\trianglelefteq\J_k(R)$.
\\

For $k\geqslant2$ and $l\geqslant k$, we will write $\J_k^l(R)$ for $\J_k(R)/\J_l(R)$.
\\

Studying $\J(R)$ as a topological group can be cumbersome, in particular for the cases where $R$ is not compact. Luckily for us, the subgroups $\J_k(\Z)$ will let us disregard this aspect. 
\\

Intuition would suggest that $\J_k(R)$ is related to 
\[\langle x+\alpha x^n|\alpha\in R,n\geqslant k\rangle.\]
The following lemma confirms some of this suspicion.

\lem \label{lem:claus} For all $k\geqslant2$, the group $\J_k(R)$ is contained in 
\[\overline{\langle x+\alpha x^n|\alpha\in R,n\geqslant k\rangle}.\]

\begin{proof}

Let $f\in\J_k(R)$. By definition, this means that 
\[f=x+\sum_{i\geqslant k}\alpha_ix^i.\]

Define the sequence $\{g_n\}_{n\in\Z^+}$ as
\[g_i:=\left\{\begin{array}{cc}
    x+\alpha_kx^k & \text{for } i\leqslant k  \\
    g_{i-1}\circ(x+(\alpha_i-\beta_{i,i-1})x^i) & \text{for } i\geqslant k+1 
\end{array}\right.,\]
where $\beta_{l,n}$ is the $l$th coefficient of $g_n$. We will prove that, for any neighborhood $V$ of $f$, there exists $n\in\Z^+$ such that 
\[g_n\in V\cap\langle x+\alpha x^n|\alpha\in R,n\geqslant k\rangle.\]

Let $V$ be a neighborhood of $f$. Since $R[[x]]$ is homeomorphic to $R^\N$ with the product topology, there exists a finite subset $I$ of $\Z_{\geqslant k}$ and a family $\{U_i\}_I$ of open subsets of $R$ such that $\alpha_i\in U_i$ for all $i\in I$, and

\[\gamma_i\in U_i, \forall i\in I\Rightarrow x+\sum_{i\geqslant k}\gamma_ix^i\in V.\]

Let $n=\max I$. Note that $\beta_{i,n}=\alpha_i$ for all $i\in\{k,\dots,n\}$. This implies that $\beta_{i,n}\in U_i$ for all $i\in I$. Thus 
\[g_n\in V\cap\langle x+\alpha x^n|\alpha\in R,n\geqslant k\rangle.\]

We conclude that $V\cap\langle x+\alpha x^n|\alpha\in R,n\geqslant k\rangle\neq\emptyset$ for any neighborhood $V$ of $f$. This is equivalent to $f\in\overline{\langle x+\alpha x^n|\alpha\in R,n\geqslant k\rangle}$. We conclude that $\J_k(R)\subseteq\overline{\langle x+\alpha x^n|\alpha\in R,n\geqslant k\rangle}$.
\end{proof}

Furthermore, if we assume $R$ to be Hausdorff, then the equality holds, giving a set of generators for these subgroups. This will let us to check many properties.

\prop \label{prop:claus} If the topology of $R$ is $\mathcal{T}_1$, then $\J_k(R)=\overline{\langle x+\alpha x^n|\alpha\in R,n\geqslant k\rangle}$.
\begin{proof}
By the latter lemma, it suffices to prove that $\overline{\langle x+\alpha x^n|\alpha\in R,n\geqslant k\rangle}\subseteq \J_k(R)$.
\\

Nota that, under the identification of $R[[x]]$ and $R^\N$ as a topological spaces, we have that $\J_k(R)$ corresponds to 
\[\pi_1^{-1}(0)\bigcap\pi_2^{-1}(1)\bigcap\left(\bigcap_{i=2}^{k-1}\pi_i^{-1}(0)\right),\]
where $\pi_i$ is the projection in the $i$-th coordinate of $R^\N$. Since $R$ is $\mathcal{T}_1$, the set $\pi_i(\alpha)$ is closed for all $i\in\Z^+$ and all $\alpha\in R$. Thus the latter set is closed. This implies that $\J_k(R)$ is closed.\\

The definition of $\J_k(R)$ implies that $x+\alpha x^n\in\J_k(R)$ for all $\alpha\in R$ and for all $n\geqslant k$. Thus
\[\langle x+\alpha x^n|\alpha\in R,n\geqslant k\rangle\subseteq\J_k(R).\]

Since $\J_k(R)$ is closed, this implies that
\[\overline{\langle x+\alpha x^n|\alpha\in R,n\geqslant k\rangle}\subseteq\J_k(R).\]
Thus concluding the proof
\end{proof}
In order to use this proposition, we will assume that $R$ is Hausdorff from here on out.
\\

We can characterize the equivalence classes in the quotients via the following proposition.

\prop Let $k$ and $l$ be positive integers such that $l>k\geqslant2$, and let $p:\J_k(R)\to\J_k^l(R)$  be the canonical projection. Let $f,g\in\J_k(R)$, with $f=x+\sum_{i\geqslant k}\alpha_ix^i$ and $g=x+\sum_{i\geqslant k}\beta_ix^i$. Then the following are equivalent:

\begin{enumerate}[$a)$]
    \item $p(f)=p(g)$.

    \item $\alpha_i=\beta_i$ for all $i\in\{k,\dots,l-1\}$.

    \item $f+x^lR[[x]]=g+x^lR[[x]]$ as elements of $R[[x]]/x^lR[[x]]$.
\end{enumerate}
\begin{proof}
    It is clear that $b)\Leftrightarrow c)$.\\
    
   We will prove that $a)\Rightarrow c)$. If $p(f)=p(g)$ then, by definition, 
   there exists $h\in\J_l(R)$ such that $f=g\circ h$. By Proposition \ref{prop:formulas}, this means that  there exists $E\in x^{k+l-1}$ such that
    \[f=g+h-x+E.\]
    Thus 
    \[f-g=h-x+E\in x^lR[[x]].\]
    This implies that $c)$ holds. \\
    
    Now, we will prove that $b)\Rightarrow a)$. If $b)$ holds, then $f-g=E\in x^lR[[x]]$. We see that
    \[g\circ f^{-1}=(f-E)(f^{-1})=f\circ f^{-1}-E(f^{-1})=x+(-E)(f^{-1})\in\J_l(R).\]
    This implies that $a)$ holds.
\end{proof}

In virtue of the latter proposition, we will write the class of $f=x+\sum_{i\geqslant k}\alpha_ix^i$ in $\J_k^l(R)$ as $f+O(x^l)$ or $x+\sum_{i=k}^{k-1}\alpha_ix^i+O(x^l)$. Furthermore, we can compute the class of $f\circ g$ with the operations of $R[[x]]/x^lR[[x]]$. This allows giving a set of generators for $\J_k^l(R)$, as stated below.

\cor For every $l \geq k$ one has 
$\J_k^l(R)=\langle x+\alpha x^n+O(x^l)|\alpha\in R,k\leqslant n \leqslant l-1\rangle.$
\begin{proof}
   Given $\Tilde{f}\in\J_k^l(R)$, let $f\in\J_k(R)$ be such that $f+O(x^l)=\Tilde{f}$. Take $\{g_n\}_{n\in\Z^+}$ as in the proof of Lemma \ref{lem:claus}. Note that $g_{l-1}$ has the same coefficients as $f$ up to the $(l-1)^{th}$ power of $x$. By the latter proposition, this implies that
   \[\Tilde{f}=f+O(x^l)=g_{l-1}+O(x^l)\in\langle x+\alpha x^n+O(x^l)|\alpha\in R,k\leqslant n \leqslant l-1\rangle.\]

   Since $\Tilde{f}$ is arbitrary, we conclude that $\J_k^l(R)\subseteq\langle x+\alpha x^n|\alpha\in R,n\geqslant k \rangle$.
\end{proof}
We now have a better understanding of the groups $\J_k(R)$, which, if a bit rudimentary, will let us move forward in their study.

\subsection{A useful lemma}
In this section we will present a result that will be central in this article.
\\

This is inspired by the work of I. K. Babenko and S. A. Bogaty\u{\i} \cite{MR2413648} and proceeds as follows: given $k\geqslant2$, we will find $c\geqslant k$ such that
$\J_c(R)\subset [\J_k(R),\J_k(R)]$.
\\

Recall that $\J_c(R)=\overline{\langle x+\alpha x^l|\alpha\in R, l\geqslant c\rangle}$. Thus, our goal translates into finding 
$c>k$ such that $x+\alpha x^l\in [\J_k(R),\J_k(R)]$ for all $\alpha\in R$, and all $l\geqslant c$. Verifying this condition directly is not easy, but we will slightly modify it to an easier-to-check yet equivalent condition.

\prop \label{prop:equi} Let $k\geqslant2$. If $c\geqslant k$, then $x+\alpha x^l\in [\J_k(R),\J_k(R)]$ for all $\alpha\in R$ and all $l\geqslant c$ if and only if for all $\alpha\in R$ and all $l\geqslant c$, there exists an element $x+\alpha x^l+\cdots$ in $[\J_k(R),\J_k(R)]$. 
\\

The proof of the latter proposition is analogous to the proof of Lemma \ref{lem:claus}. 
\\

Remember that, if $f=x+\alpha_mx^m+\alpha_{m+s}x^{m+s}+\cdots$ and $g=x+\beta_nx^n+\beta_{n+r}x^{n+r}+\cdots$ are elements of $\J(R)$ such that $m>n$ and $r>s\geqslant n$, then
\[[f,g]=x+(m-n)\alpha_m\beta_nx^{m+n-1}+C(m,n)\alpha_m\beta_n^2x^{m+2(n-1)}+\cdots.\]
The following is a consequence of this:

\prop \label{prop:imparidad} Let $k\geqslant2$.
\begin{enumerate}[$a)$]
    \item For all $\alpha\in R$ and $l\geqslant k$, there exists an element $x+\alpha x^{2l}+\cdots$ in $[\J_k(R),\J_k(R)]$.

    \item Let $l\geqslant k$ and $\alpha\in R$. Let $n\in\Z$ be such that  $n\equiv1\pmod{2}$. If there is an element $x+n\alpha x^{2l+1}+\cdots$ in $[\J_k(R),\J_k(R)]$, 
    then there is an element $x+\alpha x^{2l+1}+\cdots$ in $[\J_k(R),\J_k(R)]$.

    \item If $R$ is a domain such that $2\in R^*$, then for all $\alpha\in R$ and $l\geqslant k$, there exists an element $x+\alpha x^{2l+1}+\cdots$ in $[\J_k(R),\J_k(R)]$.
\end{enumerate}

\begin{proof}~

\begin{enumerate}[$a)$]
        \item Let $l\geqslant k$ and $\alpha\in R$. It suffices to note that
        \[[x+\alpha x^{l+1},x+x^l]=x+\alpha x^{2l}+C(l+1,l)\alpha x^{3l-1}+\cdots\]
        is such an element.

        \item By hypothesis, we know that there exists 
$m\in\Z$ such that $n=1+2m$. Furthermore 
\[[x+\alpha x^{l+2},x+x^l]=x+2\alpha x^{2l+1}+C(l+2,l)\alpha x^{3l}+\cdots.\]
Thus, we have that
\[(x+n\alpha x^{2l+1}+\cdots)\circ [x+\alpha x^{l+2},x+x^l]^{-m}=x+\alpha x^{2l+1}+\cdots.\]
Since $[\J_k(R),\J_k(R)]$ is a group, this implies that $x+\alpha x^{2l+1}+\cdots$ is contained in it.

\item Let $l\geqslant k$ and $\alpha\in R$. It suffices to note that
        \[[x+2^{-1}x^{l+2},x+\alpha x^l]=x+\alpha x^{2l+1}+C(l+2,l)2^{-1}\alpha^2x^{3l}+\cdots\]
        is such an element.
    \end{enumerate}
\end{proof}

\obs \label{obs:dos} For the case of a domain $R$ such that $2\in R^*$, this proves that 
\[\J_{2k}(R)\subset [\J_k(R),\J_k(R)].\]
We will return to this result later.
\\

For the case of a torsion-free domain $R$, we note that the elements of the form $[x+\alpha x^{l+1},x+x^l]$ and $[x+\alpha x^{l+2},x+x^l]$ will play an important role, since they will be the main tool for us to compute classes in $H_1(\J_k(R))$.
\\

Given $k\geqslant2$, the latter proposition reduces our problem to determine for which $l\geqslant k$ and $\alpha\in R$ 
there exists a product of commutators equal to some $x+n\alpha x^{2l+1}+\cdots$ with $n\equiv1\pmod{2}$.
The rest of this section is dedicated to the manufacture of such elements via the following proposition:

\prop\label{prop:cosa} Let $R$ be torsion-free domain. Given $k\geqslant2$, let $m$ and $n$ be positive integers such that $m-2>n$ and $n\geqslant k$. 
Then, for all $\alpha\in R$, there exists an element $x+C(m,n)\alpha x^{m+2(n-1)}+\cdots$ in $[\J_k(R),\J_k(R)]$.
\begin{proof}
    Let $\alpha\in R$. Note that 
    \[[x+\alpha x^m,x+x^n]=x+(m-n)\alpha x^{m+n-1}+C(m,n)\alpha x^{m+2(n-1)}+\cdots.\]
    
    If $m+n-1\equiv0\pmod{2}$, then consider the following element of $[\J_k(R),\J_k(R)]$: 
    \[[x+\alpha x^{\frac{m+n-1}{2}+1},x+x^\frac{m+n-1}{2}]^{-(m-n)}=x+(-(m-n))\alpha x^{m+n-1}+E,\]
    with $E\in x^{3\frac{m+n-1}{2}-1}R[[x]]$. We have that
    \[[x+\alpha x^m,x+x^n]\circ [x+\alpha x^{\frac{m+n-1}{2}+1},x+x^\frac{m+n-1}{2}]^{-(m-n)}=x+C(m,n)\alpha x^{m+2(n-1)}+\cdots.\]
    
    On the other hand, if $m+n-1\equiv1\pmod{2}$, then consider the following element of $[\J_k(R),\J_k(R)]$:
      \[[x+\alpha x^{\frac{m+n-2}{2}+2},x+x^{\frac{m+n-2}{2}}]^{-\frac{m-n}{2}}=x+(-(m-n))\alpha x^{m+n-1}+E,\]

      with $E\in x^{3\frac{m+n-2}{2}}R[[x]]$. We have that
       \[[x+\alpha x^m,x+x^n]\circ [x+\alpha x^{\frac{m+n-2}{2}+2},x+x^{\frac{m+n-2}{2}}]^{-\frac{m-n}{2}}= x+C(m,n)\alpha x^{m+2(n-1)}+\cdots.\]
\end{proof}

Let $\ell(m,n):=m+2(n-1)$, note that $\ell(m,n)\equiv m\pmod{2}$.
\\

Given an odd integer 
$l\geqslant 2k$ and $\alpha\in R$, the latter proposition allows us to find an element $x+\alpha x^l+\cdots$ in $[\J_k(R),\J_k(R)]$ by looking for positive integers $(m,n)$ such that $\ell(m,n)=l$ and $C(m,n)\equiv1\pmod{2}$. For this reason we will divert some paragraphs to the study of the parity of $C(m,n)$.

\prop Let $m$ and $n$ be positive integers such that $m>n$. Then $C(m,n)\equiv1\pmod{2}$ if and only if $m\equiv2\pmod{4}$ or $m\equiv3\pmod{4}$.
\begin{proof}
    Recall that 
    \[C(m,n)=\binom{m}{2}+(m-n)(m+n-1).\]
    By noting that $(m-n)(m+n-1)=m^2-m-n^2+n$, it is clear that $C(m,n)\equiv\binom{m}{2}\pmod{2}$.
\\
 
 It can be easily seen that $\binom{m}{2}\equiv1\pmod{2}$ if and only if $m\equiv2\pmod{4}$ or $m\equiv3\pmod{4}$. Thus, $C(m,n)\equiv1\pmod{2}$ if and only if $m\equiv2\pmod{4}$ or $m\equiv3\pmod{4}$.
\end{proof}

 Thanks to the latter proposition, given an odd integer $l\geqslant k$, our goal narrows down to finding a pair of positive integers $(m,n)$ such that $l=\ell(m,n),m>n+2,n\geqslant k$ and $m\equiv3\pmod{4}$.
\\

For the purpose of making the following statements easier to read, we will give a name to the pairs that verify the desired conditions:

\deff We say a pair $(m,n)\in\Z^2$ is useful if $m>n+2$ and $m\equiv1\pmod{2}$.

\deff Given $k\geqslant2$, we say the pair $(m,n)\in\Z^2$ is $k$-useful if it is useful and $n\geqslant k$. 
\\

We make the distinction between useful and $k$-useful for the sake of convenience.

\obs Let $m,n\in\Z^+$, with $m>2$. Note that $\ell(m-2,n+1)=\ell(m,n)$, and
\[C(m-2,n+1)\equiv\binom{m-2}{2}\not\equiv\binom{m}{2}\equiv C(m,n)\pmod{2}.\]
Thus, given $k\geqslant2$ and an odd integer $l\geqslant k$, our objective narrows down to finding a $k$-useful pair $(m,n)$ such that $l=\ell(m,n)$ and $(m-2,n+1)$ is too a $k$-useful pair. For this reason, we will study which are the odd numbers larger than $k$ that can be written as $\ell(m,n)$ for a $k$-useful pair, and give a lower bound to the set of such numbers.

\prop Let $l$ be an odd number. We have that: 
\begin{enumerate}[$a)$]
    \item If $l\equiv0\pmod{3}$ and $l\equiv1\pmod{4}$, then there exists $n\in\Z$ such that $n\equiv2\pmod{4}, l=\ell(n+5,n)$ and $(n+5,n)$ is a useful pair.
    
    \item If $l\equiv0\pmod{3}$ and $l\equiv3\pmod{4}$, then there exists $n\in\Z$ such that $n\equiv3\pmod{4}, l=\ell(n+8,n)$ and $(n+8,n)$ is a useful pair.
    
    \item If $l\equiv1\pmod{3}$, then there exists $n\in\Z$ such that $n\equiv1\pmod{2}, l=\ell(n+6,n)$ and $(n+6,n)$ is a useful pair.
    
    \item If $l\equiv2\pmod{3}$, then there exists $n\in\Z$ such that $n\equiv0\pmod{2}, l=\ell(n+7,n)$ and $(n+7,n)$ is a useful pair.
\end{enumerate}
\begin{proof}
~
\begin{enumerate}[$a)$]
        \item If $l\equiv0\pmod{3}$ and $l\equiv1\pmod{4}$, then by the Chinese Remainder Theorem we have that $l\equiv9\pmod{12}$. Thus, there exists $s\in\Z$ such that $l=12s+9$. Taking $n=4s+2$, we have 
        \[\ell(n+5,n)=3n+3=12s+9=l.\] 
        
        It is clear that $n\equiv2\pmod{4}$ and that $(n+5,n)$ is useful.
        
        \item If $l\equiv0\pmod{3}$  and $l\equiv3\pmod{4}$, then by the Chinese Remainder Theorem we have that  $l\equiv3\pmod{12}$. Thus, there exists $s\in\Z$ such that $l=12s+3$. Taking $n=4(s-1)+3$, we have 
       \[\ell(n+8,n)=3n+6=12s+3=l.\]
        
        It is clear that $n\equiv3\pmod{4}$ and that $(n+8,n)$ is useful.
        
        \item If $l\equiv1\pmod{3}$, then there exists $s\in\Z$ such that $l=3s+1$. Note that $l\equiv1\pmod{2}$ implies $s\equiv0\pmod{2}$. Taking $n=s-1$, we have 
         \[\ell(n+6,n)=n+6+2(n-1)=3n+3+1=3(n+1)+1=l.\]
         
        It is clear that $n\equiv1\pmod{2}$ and that$(n+6,n)$ is useful.
         
        \item If $l\equiv2\pmod{3}$, then there exists $s\in\Z$ such that $l=3s+2$. Note that $l\equiv1\pmod{2}$ implies $s\equiv1\pmod{2}$. Taking $n=s-1$, we have 
        
        \[\ell(n+7,n)=n+7+2(n-1)=3n+3+2=3(n+1)+2=l.\]  
        
        It is clear that $n\equiv0\pmod{2}$ and that $(n+7,n)$ is useful.
\end{enumerate}
\end{proof}

We are now in position to prove a first version of our desired result:

\cor Let $k\geqslant2$. If $l$ is an odd integer greater than $3k+5$, then there exists a $k$-useful pair $(m,n)$ such that $l=\ell(m,n)$ and $C(m,n)\equiv1\pmod{2}$.
\begin{proof}
 We separate by the cases of the latter proposition. We will give a proof only for case $a)$, since the proof for any other case is analogous. 
 \\
If $l\equiv0\pmod{3}$ and $l\equiv1\pmod{4}$, then there exists $n\in \Z$ such that $n\equiv2\pmod{4}$ and $l=\ell(n+5,n)$. Note that 
        \begin{align*}
            \ell(n+5,n)&\geqslant 3k+5\\
            \Leftrightarrow3(n+1)&\geqslant 3k+5\\
            \Leftrightarrow n+1&\geqslant k+\frac{5}{3}\\
            \Leftrightarrow n&\geqslant k+\frac{2}{3}\\
            \Rightarrow n&\geqslant k.
        \end{align*} 
        Thus $(n+5,n)$ is $k$-useful. Since $n+5\equiv3\pmod{4}$, we have that $C(n+5,n)\equiv1\pmod{2}$.
\end{proof}
Let $k\geqslant2$. The latter corollary, Proposition \ref{prop:cosa} and Proposition \ref{prop:imparidad} imply that, 
for all $l\geqslant3k+5$ and all $\alpha\in R$, there exists and element $x+\alpha x^l+\cdots$ in $[\J_k(R),\J_k(R)]$. By Proposition \ref{prop:equi}, we conclude that $\J_{3k+5}(R)\subset [\J_k(R),\J_k(R)]$. This bound can be improved by the same methods used to find it, but it will vary according to the class of $k$ modulo 4.\\

Define the following number:
\[c_k:=\left\{\begin{array}{cc}
    3k+1 & \text{if } k\equiv1\pmod{2} \\
    3k+2 & \text{if } k\equiv2\pmod{4} \\
    3k+4 & \text{if } k\equiv0\pmod{4}
\end{array}\right.\]
Note that  $c_k\equiv0\pmod{2}$ for all $k\geqslant2$.

\lem\label{lem:uno} If is $R$ a torsion-free domain and $k\geqslant2$, then $\J_{c_k}(R)\subset [\J_k(R),\J_k(R)]$
\begin{proof}
    We already know that $\J_{3k+5}(R)\subset [\J_k(R),\J_k(R)]$. We will prove that, for all $l\in\{c_k,\dots,3k+4\}$ and 
    all $\alpha\in R$, there exists an element $x+\alpha x^l+\cdots$ in $[\J_k(R),\J_k(R)]$. As the definition of $c_k$ suggests, we will proceed case by case.
    \begin{enumerate}[$a)$]
        \item If $k\equiv3\pmod{4}$, then we have to prove that for all $l\in\{3k+1,3k+2,3k+3,3k+4\}$ and all $\alpha\in R$, there exists an element $x+\alpha x^l+\cdots$ in $[\J_k(R),\J_k(R)]$.\\
        
        Note that $3k+1$ and $3k+3$ are even. Given that both are greater than $2k$, by Proposition \ref{prop:imparidad} we know that there are elements $x+\alpha x^{3k+1}+\cdots$ and $x+\alpha x^{3k+3}+\cdots$  in $[\J_k(\Z),\J_k(\Z)]$ for all $\alpha\in R$.\\
        
        On the other hand, note that $3k+2=\ell(k+4,k)$ and $3k+4=\ell(k+4,k+2)$. It is clear that both $(k+4,k)$ and $(k+4,k+2)$ are $k$-useful pairs. Since $k+4\equiv3\pmod{4}$, we have that
        \[C(k+4,k)\equiv1\equiv C(k+4,k+2)\pmod{2}.\]
        
       By Proposition \ref{prop:cosa} and Proposition \ref{prop:imparidad}, we conclude that for all $\alpha\in R$ there exist elements $x+\alpha x^{3k+2}+\cdots$ and $x+\alpha x^{3k+4}+\cdots$ in $[\J_k(R),\J_k(R)]$.

        \item If $k\equiv2\pmod{4}$, then we have to prove that for all $l\in\{3k+2,3k+3,3k+4\}$ and all $\alpha\in R$, there exists an element $x+\alpha x^l+\cdots$ in $[\J_k(R),\J_k(R)]$.\\
        
        Note that $3k+2$ and $3k+4$ are even. Given that both are greater than $2k$, by Proposition \ref{prop:imparidad} we know that there are elements $x+\alpha x^{3k+2}+\cdots$ and $x+\alpha x^{3k+4}+\cdots$ in $[\J_k(\Z),\J_k(\Z)]$ for all $\alpha\in R$.\\

        Moreover, note that $3k+3=\ell(k+5,k)$. It is clear that $(k+5,k)$ is a $k$-useful pair. Since $k+5\equiv3\pmod{4}$, we have that
        \[C(k+5,k)\equiv1\pmod{2}.\]
       By Proposition \ref{prop:cosa} and Proposition \ref{prop:imparidad}, we conclude that for all $\alpha\in R$ there exists an element $x+\alpha x^{3k+3}+\cdots$ in $[\J_k(R),\J_k(R)]$.

        \item If $k\equiv1\pmod{4}$, then we have to prove that for all $l\in\{3k+1,3k+2,3k+3,3k+4\}$ and all $\alpha\in R$, there exists an element $x+\alpha x^l+\cdots$ in $[\J_k(R),\J_k(R)]$.\\
        Note that $3k+1$ and $3k+3$ are even. Given that both are greater than $2k$, by Proposition \ref{prop:imparidad} we know that there are elements $x+\alpha x^{3k+1}+\cdots$ and $x+\alpha x^{3k+3}+\cdots$  in $[\J_k(\Z),\J_k(\Z)]$ for all $\alpha\in R$.\\

        Moreover, note that $3k+4=\ell(k+6,k)$. It is clear that $(k+6,k)$ is a $k$-useful pair. Since $k+6\equiv3\pmod{4}$, we have that
        \[C(k+6,k)\equiv1\pmod{2}.\]
        
       By Proposition \ref{prop:cosa} and Proposition \ref{prop:imparidad}, we conclude that for all $\alpha\in R$ there exists an element $x+\alpha x^{3k+4}+\cdots$ in $[\J_k(R),\J_k(R)]$.\\
       
        Finally, note that
        \[[x+\alpha x^{k+3},x+x^k]=x+3\alpha x^{2k+2}+C(k+3,k)\alpha x^{3k+1}+\cdots,\]
        and
        \[[x+\alpha x^{k+2},x+x^{k+1}]=x+\alpha x^{2k+2}+C(k+2,k)\alpha x^{3k+1}+\cdots.\]
        Thus, 
        \[[x+\alpha x^{k+2},x+x^{k+1}]^3\circ [x+\alpha x^{k+3},x+x^k]^{-1}=x+(3C(k+2,k)+C(k+3,k))\alpha x^{3k+1}+\cdots.\]
        Since $3C(k+2,k)+C(k+3,k)\equiv1\pmod{2}$, by Proposition \ref{prop:imparidad} we conclude that for all $\alpha\in R$ there exists an element  $x+\alpha x^{3k+1}+\cdots$ in $[\J_k(\Z),\J_k(\Z)]$.
        
        \item If $k\equiv0\pmod{4}$, then we have to prove that for all $\alpha\in R$ there exists an element $x+\alpha x^{3k+4}+\cdots$ in $[\J_k(R),\J_k(R)]$.\\
        Note that $3k+4$ is even. Given that it is greater than $2k$, by Proposition \ref{prop:imparidad} we know that for all $\alpha\in R$ there exists an element $x+\alpha x^{3k+4}+\cdots$ in $[\J_k(R),\J_k(R)]$.
        
    \end{enumerate}
\end{proof}

 With the tools at hand, the bound $c_k$ is the best we can hope among the positive integers  $l$ such that $\J_l(R)\subset [\J_k(R),\J_k(R)]$ without specifying $R$. For the case $R=\Z$, we will show that this bound is optimal, in the sense that $\J_{c_k-1}(\Z)\nsubseteq [\J_k(\Z),\J_k(\Z)]$.\\
 
The latter lemma has two main consequences, the first one of which is:
 
\obs \label{obs:diagra} Let $k\geqslant2$. Let $R$ be a torsion-free domain and $\pi:\J_k(R)\to H_1(\J_k(R))$ the natural projection. The following diagram commutes:
  \begin{center}
        \begin{tikzcd}
            \J_k(R)\arrow{rr}{\pi}\arrow{ddr}{p}&&H_1(\J_k(R))\\
            &&&\\
            &\J_k^{c_k}(R)\arrow{uur}{\Tilde{\pi}}&
        \end{tikzcd}
\end{center}
Said informally: the classes in $H_1(\J_k(R))$ can be computed using the classes of $\J_k^{c_k}(R)$. This is very useful, since $\J_k^{c_k}(R)$ is easily checked to be nilpotent.\\

The second consequence is more important. For a torsion-free domain $R$, it is not always true that $\J_{c_k}(R)= [\J_k(R),\J_k(R)]$, but the latter lemma gives a nice criteria to verify whether an element of $\J_k(R)$ is in $[\J_k(R),\J_k(R)]$ or not. With the latter remark, we have a good tool to make some headway into the question of when two elements are in the same class in $H_1(\J_k(R))$, where $R$ is a torsion-free domain.

\section{Main results}

As we discussed in the introduction, our main objective is to prove the following theorem:
\\\\
\textbf{Theorem \ref{thm:princ}}
    Let $k\geqslant2$.
    \begin{enumerate}[$a)$]
    \item  If $k\equiv 3\pmod{4}$, then \[H_1(\J_k(\Z))\cong\Z^k\bigoplus (\Z/2\Z)^\frac{k+1}{2}.\]
    
    \item If $k\equiv2\pmod{4}$, then \[H_1(\J_k(\Z))\cong\Z^k\bigoplus (\Z/2\Z)^\frac{k+2}{2}.\]
    
    \item If $k\equiv 1\pmod{4}$, then \[H_1(\J_k(\Z))\cong\Z^k\bigoplus\Z/4\Z\bigoplus (\Z/2\Z)^{\frac{k-3}{2}}.\]
    
    \item If $k\equiv 0\pmod{4}$, then \[H_1(\J_k(\Z))\cong\Z^k\bigoplus\Z/4\Z\bigoplus (\Z/2\Z)^{\frac{k-2}{2}}.\]
    \end{enumerate}

As in the previous section, we will give the results in the greatest possible generality we were able to prove. At some point we will have to restrict ourselves to $\Z$ in order to do computations.\\

For the rest of the section, we will fix $k\geqslant2$ and denote by $\pi:\J_k(R)\to H_1(\J_k(R))$ the natural projection. We will prove the main theorem via a series of intermediate results.

\subsection{Identifying the torsion of $H_1(\J_k(R))$}

As the statement of the main theorem would suggest, 
determining the torsion of $H_1(\J_k(\Z))$ is the part that will give us more problems. For this reason, we will dedicate this subsection to giving a more amenable perspective on the torsion of $H_1(\J_k(R))$, for the cases where either $R$ is a domain such that $2\in R^*$ or $R$ is a torsion-free domain. We will denote the torsion of $H_1(\J_k(R))$ by $T(H_1(\J_k(R)))$.
\\

The first result in this direction is one from the original paper \cite{MR0061610}, in which, among many other results, S. A. Jennings proved the following.
\\\\
\textbf{Proposition \ref{prop:cont}} For all $k \geqslant 2$, one has  $[\J_k(R),\J_k(R)]\subset \J_{2k}(R).$
\\

This actually yields an answer for the case of $R$ being a domain such that $2\in R^*$.

\cor \label{cor:dosss} If $R$ is a domain such that $2\in R^*$, then $[\J_k(R),\J_k(R)]=\J_{2k}(R).$
\begin{proof}
    This follows from Proposition \ref{prop:cont} and Remark \ref{obs:dos}
\end{proof}

For our other case of interest, we will have to do a bit more work. We will prove the following proposition:

\prop \label{prop:tor} If $R$ is a torsion-free domain, then $T(H_1(\J_k(R)))=\pi(\J_{2k}(R))$.
\\

In order to prove this, we will prove that $T(H_1(\J_k(R)))$ is contained in $\pi(\J_{2k}(R))$ and vice-versa. We will assume that $R$ is a torsion-free domain for the rest of the section.

\prop If $R$ is a torsion-free domain, then $T(H_1(\J_k(R)))\subseteq\pi(\J_{2k}(R))$.
\begin{proof}
    Let $\Tilde{f}\in T(H_1(\J_k(R)))$. Since $\pi$ is surjective, there exists $f\in\J_k(R)$ such that $\pi(f)=\Tilde{f}$. By definition, there exists $n\in\Z^+$ such that $\Tilde{f}^n=0$. Since 
    \[[\J_k(R),\J_k(R)]=ker(\pi),\]
     this is equivalent to $f^n\in [\J_k(R),\J_k(R)]$. By Proposition \ref{prop:cont}, this implies that $lvl(f^n)\geqslant 2k$. By Corollary \ref{cor:nvl}, we have that $lvl(f)=lvl(f^n)$. Thus $lvl(f)\geqslant2k$. This implies that $f\in\J_{2k}(R)$. Therefore $\pi(f)\in\pi(\J_{2k}(R))$. We conclude that $\pi(\J_{2k}(R))\supseteq T(H_1(\J_k(R)))$.
\end{proof}

Proving the inverse containment will require more work. We will prove that $\pi(\J_{2k}(R))$ is generated by elements of finite order.

\obs By Remark \ref{obs:diagra}, we have that
\[\pi(\J_{2k}(R))=\Tilde{\pi}(\J_{2k}^{c_k}(R)).\]
Since
\[\J_k^l(R)=\langle x+\alpha x^n+O(x^{c_k})|\alpha\in R,k\leqslant n \leqslant l-1\rangle,\]
the latter implies that
\[\Tilde{\pi}(\J_{2k}^{c_k}(R))=\langle \Tilde{\pi}(x+\alpha x^n+O(x^{c_k}))|\alpha\in R,2k\leqslant n \leqslant c_k-1\rangle.\]
Thus,
\[\pi(\J_{2k}(R))=\langle \pi(x+\alpha x^n)|\alpha\in R,2k\leqslant n\leqslant c_k-1\rangle.\]

Since $H_1(\J_k(R))$ is abelian, it will suffice to prove that $|\pi(x+\alpha x^n)|<\infty$ for all $\alpha \in R$ and for all $n\geqslant k$. For this purpose we will prove the following lemma: 

\lem If $l\geqslant2k+1$, then $\pi(x+\alpha x^l)^2=\pi(x+2\alpha x^l)$, for all $\alpha\in R$.

\begin{proof}
    Let $\alpha\in R$. We know that
    \[(x+\alpha x^l)^{~2}=x+2\alpha x^l+E,\]
    with $E\in x^{2l-1}R[[x]]$. Note that
    \[2l-1\geqslant2(2k+1)-1=4k+1.\]
    
    If $k\geqslant4$, then $4k+1>3k+4\geqslant c_k$. On the other hand, if $k\in\{2,3\}$, then it is clear that $4k+1>c_k$. Thus
    \[(x+\alpha x^l+O(x^{c_k}))^2=x+2\alpha x^l+O(x^{c_k}).\]
    This implies that
    \[\Tilde{\pi}(x+\alpha x^l+O(x^{c_k}))^2=\Tilde{\pi}(x+2\alpha x^l+O(x^{c_k})).\]
        By Remark \ref{obs:diagra}, this is equivalent to $\pi(x+\alpha x^l)^2=\pi(x+2\alpha x^l)$.
 \end{proof}

\prop\label{prop:cota} If $l\geqslant k+2$, then $\pi(x+\alpha x^{2l})=0$ and $|\pi(x+\alpha x^{2l+1})|\leqslant2$, for all $\alpha\in R$.
\begin{proof}
    Let $\alpha\in R$. Recall that, for all $l\geqslant k$, 
    \[[x+\alpha x^{l+1},x+x^l]=x+\alpha x^{2l}+C(l+1,l)\alpha x^{3l-1}+\cdots.\] 
    Thus $(x+\alpha x^{2l})\circ [x+\alpha x^{l+1},x+x^l]^{-1}\in\J_{c_k}(\Z)$ if and only if $3l-1\geqslant c_k$. Note that $l\geqslant k+2$ implies that $3l-1\geqslant3k+5$. Since $3k+5>c_k$, this implies that $\pi(f_{2l})=0$.
\\

    On the other hand, for all $l\geqslant k$,
\[(x+2\alpha x^{2l+1})\circ [x+\alpha x^{l+2},x+x^l]^{-1}\in \J_{3l}(\Z).\]
Note that $l\geqslant k+2$ implies that $3l>c_k$. Thus, $(x+2\alpha x^{2l+1})\circ [x+\alpha x^{l+2},x+x^l]^{-1}\in\J_{c_k}(\Z)$. Therefore, 
\[\pi(x+\alpha x^{2l+1})^2=\pi(x+2\alpha x^l)=0.\]
We conclude that $|\pi(x+\alpha x^{2l+1})|\leqslant2$.
\end{proof}

This determines the order of most generators and reduces their number. It also suggests some of the structure of $\pi(\J_{2k}(R))$. 
\\

At this point, we can state the following:

\prop The group $\pi(\J_{2k}(R))$ is generated by
$\{\pi(x+\alpha x^{2i+1})|k\leqslant i\leqslant\frac{c_k}{2}-1,\alpha\in R\}$.
\begin{proof}
By the latter proposition, we know that $\pi(\J_{2k}(R))$ is generated by 
\[\{\pi(x+\alpha x^{2k})|\alpha\in R\}\cup\{\pi(x+\alpha x^{2k+2})|\alpha\in R\}\cup\{\pi(x+\alpha x^{2i+1})|k\leqslant i\leqslant\tfrac{c_k}{2}-1,\alpha\in R\}.\]
We will prove that $\pi(x+\alpha_0 x^{2k})$ and $\pi(x+\alpha_0 x^{2k+2})$ are contained in $\langle\{\pi(x+\alpha x^{2i+1})|k\leqslant i\leqslant\frac{c_k}{2}-1,\alpha\in R\}\rangle$ for all $\alpha_0\in R$.
\\

    Let $\alpha_0\in R$. Consider
    \[[x+\alpha_0 x^{k+1},x+x^k]=x+\alpha_0 x^{2k}+C(k+1,k)\alpha_0 x^{3k-1}+\cdots.\]
    Note that $(x+\alpha_0 x^{2k})\circ [x+\alpha_0 x^{k+1},x+x^k]^{-1}\in \J_{3k-1}(R)$. This implies that 
    \[\pi(x+\alpha_0 x^{2k})\in \langle \{\pi(x+\alpha x^i)|3k-1\leqslant i\leqslant c_k-1,\alpha\in R\}\rangle.\]
    Since $k\geqslant2$, we have that $3k-1\geqslant2k+1$. Thus,
    \[\pi(x+\alpha_0 x^{2k})\in \langle \{\pi(x+\alpha x^i)|2k+1\leqslant i\leqslant c_k-1,\alpha\in R\}\rangle.\]
    Therefore,
    \begin{align*}
        &\langle\{\pi(x+\alpha x^{2k})|\alpha\in R\}\cup\{\pi(x+\alpha x^{2k+2})|\alpha\in R\}\cup\{\pi(x+\alpha x^{2i+1})|k\leqslant i\leqslant\tfrac{c_k}{2}-1,\alpha\in R\}\rangle\\
        =&\langle\{\pi(x+\alpha x^{2k+2})|\alpha\in R\}\cup\{\pi(x+\alpha x^{2i+1})|k\leqslant i\leqslant\tfrac{c_k}{2}-1,\alpha\in R\}\rangle.
    \end{align*}

    Analogously, note that
    \[[x+\alpha_0 x^{k+2},x+x^{k+1}]=x+\alpha_0 x^{2k+2}+C(k+2,k+1)\alpha_0 x^{3k+2}+\cdots\]
    and that $k\geqslant2$ implies $3k+2>2k+3$ . Having this, it suffices to apply the same argument as in the latter case.
\end{proof}

In order to finish this proof, we have to show that
 $\pi(x+\alpha x^{2k+1})$ and $\pi(x+\alpha x^{2k+3})$ have finite order for all $\alpha\in R$. We will distinguish 
 between the different cases given by the class of $k$ modulo $4$. This will be done this way, since it will be used in later sections.

\prop If $k\not\equiv0\pmod{4}$, then $|\pi(x+\alpha x^{2k+3})|\leqslant2$ for all $\alpha\in R$.

\begin{proof}
    Let $\alpha\in R$. Note that  
     \[[x+\alpha x^{k+3},x+x^{k+1}]=x+2\alpha x^{2k+3}+C(k+3,k+1)\alpha x^{3k+3}+\cdots.\]
     This implies that
     \[(x+2\alpha x^{2k+3})\circ [x+\alpha x^{k+3},x+x^{k+1}]^{-1}\in \J_{3k+3}(R).\]
     Since $k\not\equiv0\pmod{4}$, we know that $3k+3\geqslant c_k$. This implies that $\J_{3k+3}(R)\subseteq \J_{c_k}(R)$. Therefore, 
     \[(x+2x^{2k+3})\circ [x+\alpha x^{k+3},x+x^{k+1}]^{-1}\in \J_{c_k}(R).\]
     
     We conclude that $\pi(x+\alpha x^{2k+3})^2=\pi(x+2\alpha x^{2k+3})=0$, thus $|\pi(x+\alpha x^{2k+3})|\leqslant2$.
\end{proof}

\prop If $k\equiv0\pmod{4}$, then $|\pi(x+\alpha x^{2k+3})|$ divides $4$ for all $\alpha\in R$.

\begin{proof}
    Let $\alpha_0\in R$. Note that
    \[[x+\alpha_0 x^{k+3},x+x^{k+1}]=x+2\alpha_0x^{2k+3}+C(k+3,k+1)\alpha_0x^{3k+3}+\cdots.\]
    This implies that
    \[(x+2\alpha_0x^{2k+3})\circ [x+\alpha_0 x^{k+3},x+x^{k+1}]^{-1}\in\J_{3k+3}(R).\]
    Since $k\equiv0\pmod{4}$, we know that $\pi(\J_{3k+3}(R))=\langle\{\pi(x+\alpha x^{3k+3})|\alpha\in R\}\rangle$. Thus,
    \[\pi(x+\alpha_0x^{2k+3})^2=\pi(x+2\alpha_0x^{2k+3})\in\langle\{\pi(x+\alpha x^{3k+3})|\alpha\in R\}\rangle.\]
    This implies that there exist $\{\alpha_i\}_{i=1}^n\subset R$ and $\{\epsilon_i\}_{i=1}^n\subset \Z$ such that
    \[\pi(x+\alpha_0 x^{2k+3})^2=\prod_{i=1}^n\pi(x+\alpha_i x^{3k+3})^{\epsilon_i}.\]

    Since $3k+3\geqslant2k+4$, we have that
     \[\pi(x+\alpha_0 x^{2k+3})^4=\prod_{i=1}^n\pi(x+\alpha_i x^{3k+3})^{2\epsilon_i}=0.\]
     We conclude that $|\pi(x+\alpha_0 x^{2k+3})|$ divides $4$.
\end{proof}

\prop If $k\not\equiv3\pmod{4}$, then $|\pi(x+\alpha x^{2k+1})|$ divides $4$ for all $\alpha\in R$.

\begin{proof}
    Let $\alpha_0\in R$. Note that
    \[[x+\alpha_0 x^{k+2},x+x^k]=x+2\alpha_0x^{2k+1}+C(k+2,k)\alpha_0x^{3k}+\cdots.\]
   
    We have that
    \[(x+2\alpha_0x^{2k+1})\circ [x+\alpha_0 x^{k+2},x+x^k]^{-1}\in \J_{3k}(R).\]
    This implies that 
    \[\pi(x+\alpha_0 x^{2k+1})^2=\pi(x+2\alpha_0x^{2k+1})\in\langle\{\pi(x+\alpha x^i)|3k\leqslant i\leqslant c_k-1,\alpha\in R\}\rangle.\]
    This means that for all $i\in\{3k,\dots,c_k-1\}$, there exists $\{\alpha_{i,j}\}_{j=1}^{n_i}\subset R$ and $\{\epsilon_{i,j}\}_{j=1}^{n_i}\subset \Z$ such that
    \[\pi(x+\alpha_0 x^{2k+1})^2=\prod_{i=3k}^{c_k-1}\prod_{j=1}^{n_i}\pi(x+\alpha_{i,j} x^i)^{\epsilon_{i,j}}.\]

    We distinguish between several cases. If $k\equiv0\pmod{4}$ or $k\equiv1\pmod{4}$, then $k\geqslant4$. This implies that $3k\geqslant2k+4$. Since $H_1(\J_k(\Z))$ is abelian, we have that
    \[\pi(x+\alpha_0 x^{2k+1})^4=\left(\prod_{i=3k}^{c_k-1}\prod_{j=1}^{n_i}\pi(x+\alpha_{i,j} x^i)^{\epsilon_{i,j}}\right)^2=\prod_{i=3k}^{c_k-1}\prod_{j=1}^{n_i}\pi(x+\alpha_{i,j} x^i)^{2\epsilon_{i,j}}=0.\]
   For this case, we conclude that $|\pi(x+\alpha_0 x^{2k+1})|$ divides $4$.
    \\
    
    If $k\equiv2\pmod{4}$, then $3k\equiv0\pmod{2}$. This implies that
    \[\langle\{\pi(x+\alpha x^i)|3k\leqslant i\leqslant c_k-1,\alpha\in R\}\rangle=\langle\{\pi(x+\alpha x^i)|3k+1\leqslant i\leqslant c_k-1,\alpha\in R\}\rangle.\]
   Since $k\geqslant2$ implies that $3k+1\geqslant2k+3$, we can use an argument analogous to that of the previous case to conclude that $|\pi(x+\alpha_0 x^{2k+1})|$ divides $4$.
\end{proof}

\prop If $k\equiv3\pmod{4}$, then $|\pi(x+\alpha x^{2k+1})|\leqslant2$ for all $\alpha\in R$.

\begin{proof}
     Let $\alpha\in R$. The hypothesis implies that $c_k=3k+1$ and $C(k+2,k)\equiv0\pmod{2}$. Note that
    \[[x+\alpha x^{k+2},x+ x^k]=x+2\alpha x^{2k+1}+C(k+2,k)\alpha x^{3k}+\cdots,\]
and
     \[[x+\alpha x^{\frac{3k-1}{2}+2},x+x^{\frac{3k-1}{2}}]=x+2\alpha x^{3k}+E,\]
     with $E\in x^{3\frac{3k-1}{2}}R[[x]]$. Consider 
     \[g=[x+\alpha x^{k+2},x+x^k]^{-1}\circ [x+\alpha x^{\frac{3k-1}{2}+2},x+x^{\frac{3k-1}{2}}]^{\frac{C(k+2,k)}{2}+1}=x-2\alpha x^{2k+1}+2\alpha x^{3k}+\cdots.\]
     It is clear that $g\in [\J_k(R),\J_k(R)]$.
     \\

    On the other hand, note that
    \[(x+\alpha x^{2k+1})^2\circ g=x+2\alpha x^{3k}+\cdots.\]
    This implies that $\pi(x+\alpha x^{2k+1})^2=\pi(x+\alpha x^{3k})^2$. 
    \\
    
     Since $k\equiv3\pmod{4}$, we have that $k\geqslant3$, thus $3k\geqslant2k+3$. This implies that $\pi(x+\alpha x^{3k})^2=0$. Therefore $\pi(x+\alpha x^{2k+1})^2=0$. We conclude that $|\pi(x+\alpha x^{2k+1})|\leqslant2$.
\end{proof}

We can finally prove the following:

\prop If $R$ is a torsion-free domain, then $\pi(\J_{2k}(R))\subseteq T(H_1(\J_k(R)))$.

\begin{proof}
    Note that $\pi(\J_{2k}(R))$ is generated by elements of finite order. Since $H_1(\J_k(R))$ is abelian, this implies that every element of $\pi(\J_{2k}(R))$ has finite order. Thus $\pi(\J_{2k}(R))\subseteq T(H_1(\J_k(R)))$.
\end{proof}
With this we conclude the proof of Proposition \ref{prop:tor}.

\subsection{Approximation of the torsion of $H_1(\J_k(R))$}

Throughout this section we will assume that $R$ is a torsion-free domain. Given the results of the previous section, we will make no distinction between $T(H_1(\J_k(R)))$ and $\pi(\J_{2k}(R))$.
\\

We know that $\pi(\J_{2k}(R))$ is an abelian group, and we have a set of generators for it. We would like to find relations between these generators. For this purpose, we will reduce $\J_k^{c_k}(R)$ in order to give the structure of a finitely generated $R$-module to one of its subgroups, and show that, for some rings, the group $\pi(\J_{2k}(R))$ is a quotient of this $R$-module.
\\

We have previously made allusion to the fact that, in order to compute the class of an element $x+\sum_{i=2k}^{c_k-1}\alpha_ix^i+O(x^{c_k})$ in $H_1(\J_k(R))$, we can disregard the coefficient $\alpha_{2i}$ and reduce modulo $2$ the coefficient $\alpha_{2i+1}$, for all $i\geqslant k$. We will formalize this notion. In order to do this, we define the following number:
\[d_k:=\left\{\begin{array}{cc}
    2k+1 & \text{if } k\equiv3\pmod{4} \\
    2k+2 & \text{if } k\equiv2\pmod{4} \text{ or }k\equiv 1\pmod{4} \\
    2k+4 & \text{if } k\equiv0\pmod{4}
\end{array}\right.\]
Note that if $k\not\equiv1\pmod{4}$, then $d_k=c_k-k$.  Define the following set:

\[H_k:=\left\{x+\sum_{i=d_k}^{c_k-1}\alpha_ix^i+O(x^{c_k})\in\J_k^{c_k}(R)\bigg|\alpha_i\in 2R, \forall i\in(2\Z+1)\cap\Z_{\geqslant d_k}\right\}.\]

We will prove that this is a normal subgroup of $\J_k^{c_k}(R)$ and $\pi$ factors through the corresponding quotient.

\prop The set $H_k$ is a subgroup of $\J_k^{c_k}(R)$.

\begin{proof}
    Let $\Tilde{f},\Tilde{g}\in H_k$. By definition, we know that
     \[\Tilde{f}=x+\sum_{i=d_k}^{c_k-1}\alpha_ix^i+O(x^{c_k})\]
      and
       \[\Tilde{g}=x+\sum_{i=d_k}^{c_k-1}\beta_ix^i+O(x^{c_k}),\]
        with $\alpha_i,\beta_i\in  2R$ for all odd $i\geqslant d_k$. Let $f$ and $g$ be elements of $\J_k(R)$ such that $f+O(x^{c_k})=\Tilde{f}$ and $g+O(x^{c_k})=\Tilde{g}$. Note that
    \[f\circ g=f+g+E,\]
    with $E\in x^{2d_k-1}R[[x]]$. Since $d_k\geqslant2k+1$, we have that 
    \[2d_k-1\geqslant 2(2k+1)-1=4k+1>c_k.\]
    Thus
    \[\Tilde{f}\circ\Tilde{g}=x+\sum_{i=2k}^{c_k-1}(\alpha_i+\beta_i)x^i+O(x^{c_k}).\]
    Since $2R$ is an ideal of $R$, this implies that $\Tilde{f}\circ\Tilde{g}\in H_k$. We conclude that $H_kH_k\subset H_k$.
    \\
    
    On the other hand, let $\Tilde{f}\in H_k$. Let $f\in \J_k(R)$ be such that $f+O(x^{c_k})=\Tilde{f}$, with $f=x+\sum_{i\geqslant d_k}\alpha_ix^i$. It is clear that $f^{-1}+O(x^{c_k})=\Tilde{f}^{-1}$. We know that $lvl(f^{-1})=lvl(f)\geqslant d_k$. Thus $f^{-1}=x+\sum_{i\geqslant d_k}\beta_ix^i$. By the previous computation, we have that
    \[x+O(x^{c_k})=\Tilde{f}\circ\Tilde{f}^{-1}=x+\sum_{i=2k}^{c_k-1}(\alpha_i+\beta_i)x^i+O(x^{c_k}).\]
   Thus $\beta_i=-\alpha_i$ for all $i\in\{d_k,\dots,c_k-1\}$. Since $\alpha_i\equiv-\alpha_i\pmod{2}$ for  all  $i\geqslant d_k$, the latter equation implies that $\beta_i\in 2R$ for every odd integer $i$ in $\{d_k,\dots,c_k-1\}$. The fact that $\tilde{f}\in H_k$ follows. We conclude that $(H_k)^{-1}\in H_k$.
    \\
    
    Finally, since $0\in 2R$, we have that $x+O(x^{c_k})\in H_k$.
\end{proof}

\prop The subgroup $H_k$ is normal in $\J_k^{c_k}(R)$.
\begin{proof}
    Let $\Tilde{f}\in\J_k^{c_k}(R)$ and $\Tilde{g}\in H_k$. Let $f$ and $g$ be elements in $\J_k(R)$ such that $f+O(x^{c_k})=\Tilde{f}$ and $g+O(x^{c_k})=\Tilde{g}$. Note that
    \[f\circ g\circ f^{-1}=[f,g]\circ g.\]
    
    By definition, we have that $lvl(g)\geqslant d_k$. Furthermore, we know that 
    \[lvl([f,g])\geqslant lvl(f)+lvl(g)-1.\]
    
    On the one hand, if $lvl(f)\geqslant k+1$ or $lvl(g)\geqslant d_k+1$, then
     \[lvl([f,g])\geqslant k+d_k=c_k.\]
     Thus $[\Tilde{f},\Tilde{g}]=x+O(x^{c_k})$. This implies that $\Tilde{f}\circ\Tilde{g}\circ\Tilde{f}^{-1}=\Tilde{g}\in H_k$.
     \\
     
     On the other hand, if $k\not\equiv1\pmod{4}, lvl(f)=k$ and $lvl(g)=d_k$, then $[\Tilde{f},\Tilde{g}]=x+(k-d_k)x^{c_k-1}+O(x^{c_k})$. This implies that
     \[\Tilde{f}\circ\Tilde{g}\circ\Tilde{f}^{-1}=(x+(k-d_k)x^{c_k-1}+O(x^{c_k}))\circ g=g-(d_k-k)x^{c_k-1}+O(x^{c_k}).\]
      Note that $d_k-k=c_k-2k$. Since $c_k\equiv0\pmod{2}$, we have that $d_k-k\equiv0\pmod{2}$. This implies that $\Tilde{f}\circ\Tilde{g}\circ\Tilde{f}^{-1}\in H_k$.
     \\
     
     Finally, if $k\equiv1\pmod{4}, lvl(f)=k$ and $lvl(g)=d_k$, then we can apply the same argument as in the first case.

     We conclude that $H_k$ is normal.
\end{proof}
Define $\Gamma_k:=\J_k^{c_k}(R)/H_k$.

\prop The subgroup $H_k$ is contained $ker(\Tilde{\pi})$.

\begin{proof}
    Let 
    \[E:=\{x+\alpha x^i+O(x^{c_k})|\alpha\in R, d_k\leqslant i\leqslant c_k-2,\text{ and $i$ is even}\},\]
    and
    \[O:=\{(x+\alpha x^i+O(x^{c_k}))^2|\alpha\in R, d_k\leqslant i\leqslant c_k-1,\text{ and $i$ is odd}\}.\]    
    Note that 
    \[H_k=\langle E\cup O\rangle.\]
    
    We have that
    \begin{align*}
        E
        =\{x+\alpha x^i+O(x^{c_k})|\alpha\in R, d_k\leqslant i\leqslant 2k+2,\text{ and $i$ is even}\}\cup\{x+\alpha x^{2i}+O(x^{c_k})|\alpha\in R, k+2\leqslant i\leqslant \tfrac{c_k}{2}-1\},
    \end{align*}
    and
    \begin{align*}
        &O\\
        =&\{(x+\alpha x^i+O(x^{c_k}))^2|\alpha\in R, d_k\leqslant i\leqslant 2k+3,\text{ and $i$ is odd}\}\cup\{(x+\alpha x^{2i+1}+O(x^{c_k}))^2|\alpha\in R, k+2\leqslant i\leqslant \tfrac{c_k}{2}-1\}.
    \end{align*}

    Since $\Tilde{\pi}\circ p=\pi$, we have that $\pi(x+\alpha x^i)=\Tilde{\pi}(x+\alpha x^i+O(x^{c_k}))$ for all $\alpha\in R$ and for all $i\geqslant k$. By Proposition \ref{prop:cota}, this implies that 
    \[\{x+\alpha x^{2i}+O(x^{c_k})|\alpha\in R, k+2\leqslant i\leqslant \tfrac{c_k}{2}-1\}\subset ker(\Tilde{\pi}).\]
    and
     \[\{(x+\alpha x^{2i+1}+O(x^{c_k}))^2|\alpha\in R, k+2\leqslant i\leqslant \tfrac{c_k}{2}-1\}\subset ker(\Tilde{\pi}).\]

    Thus, it suffices to prove that $\Tilde{\pi}(x+\alpha x^i+O(x^{c_k}))=0$ for all $\alpha\in R$ and for every even integer $i$ in $\{d_k,\dots, 2k+2\}$, 
    and that $\Tilde{\pi}(x+\alpha x^i+O(x^{c_k}))^2=0$ for all $\alpha\in R$
    and for every odd integer $i$ in $\{d_k,\dots, 2k+3\}$. For this purpose we will take cases depending on the class modulo $4$ of $k$.

     \begin{enumerate}[$a)$]
         \item If $k\equiv3\pmod{4}$, then $d_k=2k+1$. We know that 
         \[\Tilde{\pi}(x+\alpha x^{2k+1}+O(x^{c_k}))^2=0=\Tilde{\pi}(x+\alpha x^{2k+3}+O(x^{c_k}))^2,\]
         for all $\alpha\in R$. Thus, it will suffice to prove that $\Tilde{\pi}(x+\alpha x^{2k+2})=0$ for all $\alpha\in R$.
         \\
         
         Let $\alpha\in R$. We have that
         \[[x+\alpha x^{k+2},x+x^{k+1}]=x+\alpha x^{2k+2}+C(k+2,k+1)\alpha x^{3k+3}+\cdots.\]
         Note that $(x+\alpha x^{2k+2})\circ [x+\alpha x^{k+2},x+x^{k+1}]^{-1}\in \J_{3k+3}(R)\subset \J_{c_k}(R)$. This implies that 
         \[\Tilde{\pi}(x+\alpha x^{2k+2}+O(x^{c_k}))=\pi(x+\alpha x^{2k+2})=0.\]

         \item If $k\equiv2\pmod{4}$ or $k\equiv1\pmod{4}$, then $d_k=2k+2$. We know that $\Tilde{\pi}(x+\alpha x^{2k+3}+O(x^{c_k}))^2=0$ for all $\alpha\in R$. The proof that $\Tilde{\pi}(x+\alpha x^{2k+2})=0$ for all $\alpha\in R$ is analogous to the previous case.

         \item If $k\equiv0\pmod{4}$, then $d_k=2k+4$ and we are already done.
         
     \end{enumerate}
Since $E\cup O\subset ker(\Tilde{\pi})$, we conclude  that $H_k\subset ker(\Tilde{\pi})$.
    
\end{proof}

We can now apply the factor theorem to obtain a unique morphism $\rho:\Gamma_k\to H_1(\J_k(R))$ such that the following diagram commutes:
    \begin{center}
        \begin{tikzcd}
            \J_k^{c_k}(R)\arrow{rr}{\Tilde{\pi}}\arrow{ddr}{p'}&&H_1(\J_k(R))\\
            &&&\\
            &\Gamma_k\arrow{uur}{\rho}&
        \end{tikzcd}
\end{center}
Moreover, we can compose this diagram with the diagram from Remark \ref{obs:diagra} to obtain the following commutative diagram:

    \begin{center}
        \begin{tikzcd}
            \J_k(R)\arrow{rr}{\pi}\arrow{ddr}{p'\circ p}&&H_1(\J_k(R))\\
            &&&\\
            &\Gamma_k\arrow{uur}{\rho}&
        \end{tikzcd}
\end{center}

From any of the previous diagrams, it is clear that $\pi(\J_{2k}(R))=¸\rho\circ p'(\J_{2k}^{c_k}(R))$. We endow $p'(\J_{2k}^{c_k}(R))$ with the structure of an $R$-module. 

\prop Let $f$ and $g$ be elements in $\Gamma_k$, with $f=p'(x+\sum_{i=2k}^{c_k-1}\alpha_ix^i+O(x^{c_k}))$ and $g=p'(x+\sum_{i=2k}^{c_k-1}\beta_ix^i+O(x^{c_k}))$. Then

\[f\circ g=p'\left(x+\sum_{i=2k}^{c_k-1}(\alpha_i+\beta_i)x^i+O(x^{c_k}))\right).\]

\begin{proof}
    From Proposition \ref{prop:formulas} we have that 

   \[\left(x+\sum_{i=2k}^{c_k-1}\alpha_ix^i\right)\circ \left(x+\sum_{i=2k}^{c_k-1}\beta_ix^i\right)=x+\sum_{i=2k}^{c_k-1}(\alpha_i+\beta_i)x^i+2k\alpha_{2k}\beta_{2k}x^{4k-1}+E,\]
   with $E\in x^{6k-2}R[[x]]$. Thus
   \begin{align*}
       f\circ g=&p'\left(x+\sum_{i=2k}^{c_k-1}\alpha_ix^i+O(x^{c_k})\right)\circ p'\left(x+\sum_{i=2k}^{c_k-1}\beta_ix^i+O(x^{c_k})\right)\\
       =&p'\left(x+\sum_{i=2k}^{c_k-1}(\alpha_i+\beta_i)x^i+2k\alpha_{2k}\beta_{2k}x^{4k-1}+O(x^{c_k})\right).
   \end{align*}
   Note that $x+2\alpha_{2k}\beta_{2k}x^{4k-1}+O(x^{c_k})\in H_k$. Moreover 
   \[x+\sum_{i=2k}^{c_k-1}(\alpha_i+\beta_i)x^i+2k\alpha_{2k}\beta_{2k}x^{4k-1}+O(x^{c_k})=\left(x+\sum_{i=2k}^{c_k-1}(\alpha_i+\beta_i)x^i+O(x^{c_k})\right)\circ(x+2\alpha_{2k}\beta_{2k}x^{4k-1}+O(x^{c_k})).\]
   This implies the equality.
\end{proof}

From the previous proposition, it is clear that
\begin{align*}
    R\times p'(\J_{2k}^{c_k}(R))&\to p'(\J_{2k}^{c_k}(R))\\
    \left(r,p'\left(x+\sum_{i=2k}^{c_k-1}\alpha_ix^i+O(x^{c_k})\right)\right)&\mapsto p'\left(x+\sum_{i=2k}^{c_k-1}r\alpha_ix^i+O(x^{c_k})\right)
\end{align*}
defines the structure of an $R$-module on $p'(\J_{2k}^{c_k}(R))$. In order to translate this  $R$-module structure onto $\pi(\J_{2k}(R))$, we need to prove that $[\Gamma_k,\Gamma_k]$ is an $R$-submodule of $p'(\J_{2k}^{c_k}(R))$. This is not necessarily the case for any ring, but it will hold for the case of $R=\Z$. If such a condition is met, we can use the following proposition:

\prop \label{prop:clasi} Let $R$ be a torsion-free PID. If $[\Gamma_k,\Gamma_k]$ is an $R$-submodule of $p'(\J_{2k}^{c_k}(R))$, then 
\[T(\J_k(R))\cong \bigoplus_{i=1}^{n_k}R/(p_i),\]
where $n_k\leqslant\frac{c_k-2k}{2}$ and $p_i|4$ for all $i\in\{1,\dots,n_k\}$.

\begin{proof}
    By hypothesis, we have that 
    \[p'(\J_{2k}^{c_k}(R))/[\Gamma_k,\Gamma_k]\cong T(\J_k(R))\]    
    is an $R$-module. Since 
    \[\{\pi(x+\alpha x^{2i+1})|\alpha\in R,k\leqslant i \leqslant \tfrac{c_k}{2}-1\}\]
    generate $T(H_1(\J_k(R)))$ as a group, it follows that
    \[\{\pi(x+x^{2i+1})|k\leqslant i \leqslant \tfrac{c_k}{2}-1\}\]
    generates $T(H_1(\J_k(R)))$ as an $R$-module. Note that 
    \[\#\{\pi(x+x^{2i+1})|k\leqslant i \leqslant \tfrac{c_k}{2}-1\}=\frac{c_k-2k}{2}.\]

    Since $4\pi(x+x^{2i+1})=0$ for all $i \in \{k,\dots,\frac{c_k}{2}-1\}$, the structure theorem for finitely generated modules over PID implies that
    \[T(H_1(\J_k(R)))\cong \bigoplus_{i=1}^{n_k}R/(p_i),\]
where $n_k\leqslant\frac{c_k-2k}{2}$ and $p_i|4$ for all $i\in\{1,\dots,n_k\}$.
\end{proof}

\subsection{Computing the torsion of $H_1(\J_k(\Z))$}

We will finally focus on our case of interest, thus we fix $R=\Z$. As the table of contents suggests, we will proceed by cases. These will be, roughly speaking, the general case, and the three non-general cases. Regardless of this, the method used will be, roughly speaking again, to describe $[\Gamma_k,\Gamma_k]$ as explicitly as possible, then figure out which generators are related, and finally determine its invariant factors as a $\Z$-module.

\subsubsection{The case $k\geqslant5$.}

From now on, we fix $k\geqslant5$. By Corollary \ref{cor:conmu} we have the following:

\prop Let $f$ and $g$ be elements in $\Gamma_k$, with $f=p'(x+\sum_{i=k}^{c_k-1}\alpha_ix^i+O(x^{c_k}))$ and $g=p'(x+\sum_{i=k}^{c_k-1}\beta_ix^i+O(x^{c_k}))$.
    \begin{enumerate}[$a)$]
        \item If $k\equiv3\pmod{4}$, then
        \[[f,g]=p'(x+(\alpha_{k+1}\beta_k-\alpha_k\beta_{k+1})x^{2k}+(\alpha_{k+1}\beta_k-\alpha_k\beta_{k+1})x^{3k}+O(x^{3k+1})).\]

        \item If $k\equiv2\pmod{4}$, then
        \[[f,g]=p'(x+(\alpha_{k+1}\beta_k-\alpha_k\beta_{k+1})x^{2k}+2(\alpha_{k+2}\beta_k-\alpha_k\beta_{k+2})x^{2k+1}+(\alpha_{k+1}\beta_k-\alpha_k\beta_{k+1})x^{3k-1}+O(x^{3k+2})).\]

        \item If $k\equiv1\pmod{4}$, then
        \begin{align*}[f,g]=p'(
            x+(\alpha_{k+1}\beta_k-\alpha_k\beta_{k+1})x^{2k}+2(\alpha_{k+2}\beta_k-\alpha_k\beta_{k+2})x^{2k+1}+(\alpha_{k+2}\beta_k-\alpha_k\beta_{k+2})x^{3k}+O(x^{3k+1})).
        \end{align*}

        \item If $k\equiv0\pmod{4}$, then
        \begin{align*}
            [f,g]=p'(x&+(\alpha_{k+1}\beta_k-\alpha_k\beta_{k+1})x^{2k}+2(\alpha_{k+2}\beta_k-\alpha_k\beta_{k+2})x^{2k+1}\\
            +&(3(\alpha_{k+3}\beta_k-\alpha_k\beta_{k+3})+(\alpha_{k+2}\beta_{k+1}-\alpha_{k+1}\beta_{k+2}))x^{2k+2}\\
            +&(4(\alpha_{k+4}\beta_k-\alpha_k\beta_{k+4})+2(\alpha_{k+3}\beta_{k+1}-\alpha_{k+1}\beta_{k+3}))x^{2k+3}\\
            +&(\alpha_{k+3}\beta_k+\beta_{k+3}\alpha_k)x^{3k+1}+(\alpha_{k+3}\beta_{k+1}-\beta_{k+3}\alpha_{k+1})x^{3k+3}+O(x^{3k+4})).
        \end{align*}       
    \end{enumerate}

We use these formulas to describe $[\Gamma_k,\Gamma_k]$ through sets of generators.

\cor ~
\begin{enumerate}[$a)$]
    \item If $k\equiv3\pmod{4}$, then $[\Gamma_k,\Gamma_k]=\langle p'(x+x^{2k}+x^{3k}+O(x^{3k+1}))\rangle$.

    \item If $k\equiv2\pmod{4}$, then $[\Gamma_k,\Gamma_k]=\langle p'(x+x^{2k}+x^{3k-1}+O(x^{3k+2})), p'(x+2x^{2k+1}+O(x^{3k+2}))\rangle$.

    \item If $k\equiv1\pmod{4}$, then $[\Gamma_k,\Gamma_k]=\langle p'(x+x^{2k}+O(x^{3k+1})), p'(x+2x^{2k+1}+x^{3k}+O(x^{3k+1}))\rangle.$
    
    \item If $k\equiv0\pmod{4}$, then 
    \begin{align*}
        [\Gamma_k,\Gamma_k]=&\langle p'(x+x^{2k}+O(x^{3k+4})), p'(x+2x^{2k+1}+O(x^{3k+4})), p'(x+x^{2k+2}+O(x^{3k+4})), \\
        &p'(x+x^{3k+1}+O(x^{3k+4})),p'(x+4x^{2k+3}+O(x^{3k+4})) p'(x+2x^{2k+3}+x^{3k+3}+O(x^{3k+4}))\rangle
    \end{align*}
\end{enumerate}
\begin{proof}

Most of the generators are easily seen to be contained in their respective cases, and the equality comes from factoring the respective generators with their respective powers. 
The only case were it is not obvious that this is the case is for $k\equiv0\pmod{4}$ and the generator $p'(x+x^{3k+1}+O(x^{3k+4}))$. For this case, it suffices to note that
\begin{align*}
    &[ p'(x+x^{k+3}+O(x^{3k+1})), p'(x+x^k+O(x^{3k+1}))]\circ [ p'(x+x^{k+2}+O(x^{3k+1})), p'(x+x^{k+1}+O(x^{3k+1}))]^{-3}\\
    =& p'(x+x^{3k+1}+O(x^{3k+1})).
\end{align*}
This concludes the proof.
\end{proof}

This allow us to determine the relations of the generators of $\pi(\J_{2k}(\Z))$.

\prop~
\begin{enumerate}[$a)$]
    \item If $k\equiv3\pmod{4}$, then, for $i\in\{k,\dots,\frac{3k-1}{2}\}$ we have that 
    \[\pi(x+x^{2i+1})\not\in\left\langle\{\pi(x+x^{2j+1})\}_{j=k}^\frac{3k-1}{2}\setminus\{\pi(x+x^{2i+1})\}\right\rangle.\]

    \item If $k\equiv2\pmod{4}$, then $|\pi(x+x^{2k+1})|=2$, and for $i\in\{k,\dots,\frac{3k}{2}\}$ we have that 
    \[\pi(x+x^{2i+1})\not\in\left\langle\{\pi(x+x^{2j+1})\}_{j=k}^\frac{3k}{2}\setminus\{\pi(x+x^{2i+1})\}\right\rangle.\]

    \item If $k\equiv1\pmod{4}$, then $\pi(x+x^{2k+1})^2=\pi(x+x^{3k})$, and for $i\in\{k,\dots,\frac{3k-3}{2}\}$ we have that 
    \[\pi(x+x^{2i+1})\not\in\left\langle\{\pi(x+x^{2j+1})\}_{j=k}^\frac{3k-3}{2}\setminus\{\pi(x+x^{2i+1})\}\right\rangle.\]

    \item If $k\equiv0\pmod{4}$, then $\pi(x+x^{2k+3})^2=\pi(x+x^{3k+3}),|\pi(x+x^{2k+1})|=2,\pi(x+x^{3k+1})=0$, and for $i\in\{k,\dots,\frac{3k-2}{2}\}$ we have that
    \[\pi(x+x^{2i+1})\not\in\left\langle\{\pi(x+x^{2j+1})\}_{j=k}^\frac{3k-2}{2}\setminus\{\pi(x+x^{2i+1})\}\right\rangle.\]
\end{enumerate}
\begin{proof}
We will only give the proof for the case of $k\equiv3\pmod{4}$, since for any other case the statement is either direct from the previous corollary or analogous to this case.\\

    Recall that, for $f\in\J_k(\Z)$, we 
have that $\pi(f)=\rho\circ p'(f+O(x^{c_k}))$.\\

Let $i\in\{k,\dots,\frac{3k-1}{2}\}$. Note that
        \[\pi(x+x^{2i+1})\in\left\langle\{\pi(x+x^{2j+1})\}_{j=k}^\frac{3k-1}{2}\setminus\{\pi(x+x^{2i+1})\}\right\rangle\]        
         holds if and only if for all $j\in\{k,\dots,\frac{3k-1}{2}\}\setminus\{i\}$ there exists $\epsilon_j\in\{0,1\}$ such that  
        \[\pi(x+x^{2i+1})=\prod_{\substack{k\leqslant j\leqslant\frac{3k-1}{2}\\j\neq i}}\pi(x+x^{2j+1})^{\epsilon_j}=\pi\left(x+\sum_{\substack{k\leqslant j\leqslant\frac{3k-1}{2}\\j\neq i}}\epsilon_jx^{2j+1}\right).\]        
        This equation holds if and  only  if there exists $g\in [\Gamma_k,\Gamma_k]$ such that        
        \[ p'(x+x^{2i+1}+O(x^{3k+1}))g=p'\left(x+\sum_{\substack{k\leqslant j\leqslant\frac{3k-1}{2}\\j\neq i}}\epsilon_jx^{2j+1}+O(x^{3k+1})\right).\]
        
        By the latter corollary, there exists $\alpha\in\Z$ such that         
        \[g= p'(x+x^{2k}+x^{3k}+O(x^{3k+1}))^\alpha= p'(x+\alpha x^{2k}+\alpha x^{3k}+O(x^{3k+1})).\]        
        Thus, the equation holds if and only if there exists $\alpha\in\Z$ such that         
         \[ p'(x+\alpha x^{2k}+x^{2i+1}+\alpha x^{3k}+O(x^{3k+1}))=p'\left(x+\sum_{\substack{k\leqslant j\leqslant\frac{3k-1}{2}\\j\neq i}}\epsilon_jx^{2j+1}+O(x^{3k+1})\right).\]
        
        By definition of  $p'$, the latter equation holds if and only if there exists $h=\sum_{l=2k+1}^{3k}\beta_lx^l$ such that $l\equiv1\pmod{2}$ implies  $\beta_l\equiv0\pmod{2}$ and        
       \[ x+\alpha x^{2k}+x^{2i+1}+\alpha x^{3k}+h+O(x^{3k+1})=x+\sum_{\substack{k\leqslant j\leqslant\frac{3k-1}{2}\\j\neq i}}\epsilon_jx^{2j+1}+O(x^{3k+1}).\]        
        This implies that $\alpha=0$. Thus,
        \[ x+x^{2i+1}+h+O(x^{3k+1})=x+\sum_{\substack{k\leqslant j\leqslant\frac{3k-1}{2}\\j\neq i}}\epsilon_jx^{2j+1}+O(x^{3k+1}).\]       
       Note that the $(2i+1)$th coefficient on the left side of the equation is odd. On the  other hand the $(2i+1)$th coefficient on the right side of the equation is zero, which is even. This is a contradiction.
\end{proof}

Note that the latter proposition implies that, except for the elements explicitly shown, none of the generators is in the class of $0$.\\

We can finally compute  the torsion of $H_1(\J_k(\Z))$.

\thm \label{thm:tors} Let $k\geqslant5$.
\begin{enumerate}[$a)$]
    \item If $k\equiv3\pmod{4}$, then $T(H_1(\J_k(\Z)))\cong(\Z/2 \Z)^\frac{k+1}{2}$.

    \item If $k\equiv2\pmod{4}$, then $T(H_1(\J_k(\Z)))\cong(\Z/2 \Z)^{\frac{k+2}{2}}$.

    \item If $k\equiv1\pmod{4}$, then $T(H_1(\J_k(\Z)))\cong\Z/4\Z\bigoplus(\Z/2 \Z)^\frac{k-3}{2}$.

    \item if $k\equiv0\pmod{4}$, then $T(H_1(\J_k(\Z)))\cong\Z/4\Z\bigoplus(\Z/2 \Z)^{\frac{k-2}{2}}$.
\end{enumerate}
\begin{proof} 
Recall that $T(H_1(\J_k(\Z)))=\pi(\J_{2k}(\Z))$. We prove only the case for $k\equiv3\pmod{4}$, since the other cases are analogous.\\

Note that $\pi(\J_{2k}(\Z))=\langle \pi(x+x^{2i+1})|k\leqslant i\leqslant \frac{3k-1}{2}\rangle$, and this set of generators is minimal by the latter proposition. By Proposition \ref{prop:clasi}, we have 
\[T(H_1(\J_k(\Z)))\cong \bigoplus_{i=1}^{\frac{k+1}{2}}\Z/p_i\Z,\]
with $p_i\in\{1,2,4\}$ for all $i\in\{1,\dots,\frac{k+1}{2}\}$. Note that $p_i\not\in  \{1,4\}$, otherwise we would have a null generator or a generator whose order is greater than $2$, respectively. We conclude that
\[T(H_1(\J_k(\Z)))\cong \bigoplus_{i=1}^{\frac{k+1}{2}}\Z/2\Z=(\Z/2\Z)^\frac{k+1}{2}.\]
\end{proof}

\subsubsection{The case $k\in\{2,4\}$}
By analogous methods to those of the general case we can prove the following:

\prop Let $f,g\in\Gamma_2$, with $f=p'(x+\sum_{i=2}^7\alpha_ix^i+O(x^8))$ and $g=p'(x+\sum_{i=2}^7\beta_ix^i+O(x^8))$. Then
\begin{align*}
    f\circ g=p'(x&+(\alpha_2+\beta_2)x^2+(\alpha_3+\beta_3+2\alpha_2\beta_2)x^3+(\alpha_4+\beta_4+3\alpha_3\beta_2+2\alpha_2\beta_3+\alpha_2\beta_2^{~2})x^4\\
    +&(\alpha_5+\beta_5+4\alpha_4\beta_2+3\alpha_3\beta_3+2\alpha_2\beta_4+3\alpha_3\beta_2^{~2}+2\alpha_2\beta_2\beta_3)x^5\\
    +&(\alpha_7+\beta_7+\alpha_5\beta_3+\alpha_3\beta_3+\alpha_3\beta_5+\alpha_3\beta_2\beta_3)x^7+O(x^8)).
\end{align*}

The latter formula can be obtained by using Lemma 2.2 of \cite{MR2413648} and making the respective reductions in $\Gamma_2$.

\cor Let $f\in\Gamma_2$, with $f=p'(x+\sum_{i=2}^7\alpha_ix^i+O(x^8))$. Then
\begin{align*}
    f^{-1}=p'(x&+(-\alpha_2)x^2+(-\alpha_3+2\alpha_2^{~2})x^3+(-\alpha_4+5\alpha_3\alpha_2-5\alpha_2^{~3})x^4\\
    +&(-\alpha_5+6\alpha_4\alpha_2+14\alpha_2^{~4}+3\alpha_3^{~2}-21\alpha_3\alpha_2^{~2})x^5+\alpha_7x^7+O(x^8)).
\end{align*}

\cor Let $f,g\in\Gamma_2$, with $f=p'(x+\sum_{i=2}^7\alpha_ix^i+O(x^8))$ and $g=p'(x+\sum_{i=2}^7\beta_ix^i+O(x^8))$. Then
\begin{align*}
    [f,g]=p'(x&+(\alpha_2\beta_2(\beta_2-\alpha_2)+(\alpha_3\beta_2-\alpha_2\beta_3))x^4\\
    +&(2\alpha_2\beta_2(2(\alpha_2^{~2}-\beta_2^{~2})+3(\beta_3-\alpha_3))+2(\alpha_4\beta_2-\alpha_2\beta_4)-(\alpha_3\beta_2^{~2}-\beta_3\alpha_2^{~2})x^5\\
    +&(\alpha_3\beta_2-\alpha_2\beta_3)x^7+O(x^8)).
\end{align*}

We will not detect a set of generators for $[\Gamma_2,\Gamma_2]$ through this formula, but it will allow us to describe any element of this subgroup. We will make use of this in the 
proof of the next proposition.

\prop~ The following holds: 
\begin{enumerate}[$a)$]
    \item $|\pi(x+x^5)|=2$.

    \item $\pi(x+x^7)\not\in\langle\pi(x+x^5)\rangle$.
\end{enumerate}
\begin{proof}
    Note that if $f\in [\Gamma_2,\Gamma_2]$ then, by definition, there exists $n\in\Z^+$ such that $f=\prod_{i=1}^n[f_i,g_i]$, with $f_i,g_i\in\Gamma_2$ for all $i\in\{1,\dots,n\}$. By the latter corollary we have that, if $f=p'(x+\sum_{j=2}^7\alpha_{i,j}x^j+O(x^8))$ and $g=p'(x+\sum_{j=2}^7\beta_{i,j}x^j+O(x^8))$, then
    \begin{align*}
    f=\prod_{i=1}^np'(x&+(\alpha_{i,2}\beta_{i,2}(\beta_{i,2}-\alpha_{i,2})+(\alpha_{i,3}\beta_{i,2}-\alpha_{i,2}\beta_{i,3}))x^4\\
    +&(2\alpha_{i,2}\beta_{i,2}(2(\alpha_{i,2}^{~~~2}-\beta_{i,2}^{~~~2})+3(\beta_{i,3}-\alpha_{i,3}))+2(\alpha_{i,4}\beta_{i,2}-\alpha_{i,2}\beta_{i,4})-(\alpha_{i,3}\beta_{i,2}^{~~~2}-\beta_{i,3}\alpha_{i,2}^{~~~2})x^5\\
    +&(\alpha_{i,3}\beta_{i,2}-\alpha_{i,2}\beta_{i,3})x^7+O(x^8))\\
    =p'\Bigg(x&+\left(\sum_{i=1}^n(\alpha_{i,2}\beta_{i,2}(\beta_{i,2}-\alpha_{i,2})+(\alpha_{i,3}\beta_{i,2}-\alpha_{i,2}\beta_{i,3}))\right)x^4\\
    +&\left(\sum_{i=1}^n(2\alpha_{i,2}\beta_{i,2}(2(\alpha_{i,2}^{~~~2}-\beta_{i,2}^{~~~2})+3(\beta_{i,3}-\alpha_{i,3}))+2(\alpha_{i,4}\beta_{i,2}-\alpha_{i,2}\beta_{i,4})-(\alpha_{i,3}\beta_{i,2}^{~~~2}-\beta_{i,3}\alpha_{i,2}^{~~~2})\right)x^5\\
    +&\left(\sum_{i=1}^n(\alpha_{i,3}\beta_{i,2}-\alpha_{i,2}\beta_{i,3})\right)x^7+O(x^8)\Bigg).
\end{align*}
We will make use of this formula later in this proof.
\begin{enumerate}[$a)$]
    \item Note that  $[p'(x+x^4+O(x^8)),p'(x+x^2+O(x^8))]=p'(x+2x^5+O(x^8))$. Thus,
    \[\pi(x+x^5)^2=\pi(x+2x^5)=\rho\circ p'(x+2x^5+O(x^8))=0.\]
    This implies that $|\pi(x+x^5)|\leqslant2$. Furthermore, note that $\pi(x+x^5)=0$ if and only if $p'(x+x^5+O(x^8))\in [\Gamma_2,\Gamma_2]$. It follows that
    \begin{align*}
    p'(x+x^5+O(x^8))=&p'\Bigg(x+\left(\sum_{i=1}^n(\alpha_{i,2}\beta_{i,2}(\beta_{i,2}-\alpha_{i,2})+(\alpha_{i,3}\beta_{i,2}-\alpha_{i,2}\beta_{i,3}))\right)x^4\\
    +&\Bigg(\sum_{i=1}^n(2\alpha_{i,2}\beta_{i,2}(2(\alpha_{i,2}^{~~~2}-\beta_{i,2}^{~~~2})+3(\beta_{i,3}-\alpha_{i,3}))\\
    +&2(\alpha_{i,4}\beta_{i,2}-\alpha_{i,2}\beta_{i,4})-(\alpha_{i,3}\beta_{i,2}^{~~~2}-\beta_{i,3}\alpha_{i,2}^{~~~2})\Bigg)x^5\\
    +&\left(\sum_{i=1}^n(\alpha_{i,3}\beta_{i,2}-\alpha_{i,2}\beta_{i,3})\right)x^7+O(x^8)\Bigg).
\end{align*}

On the one hand, note that 
\[\sum_{i=1}^n(2\alpha_{i,2}\beta_{i,2}(2(\alpha_{i,2}^{~~~2}-\beta_{i,2}^{~~~2})+3(\beta_{i,3}-\alpha_{i,3}))+2(\alpha_{i,4}\beta_{i,2}-\alpha_{i,2}\beta_{i,4})-(\alpha_{i,3}\beta_{i,2}^{~~~2}-\beta_{i,3}\alpha_{i,2}^{~~~2})=1.\]

On the other hand, note that $\sum_{i=1}^n(\alpha_{i,3}\beta_{i,2}-\alpha_{i,2}\beta_{i,3})\equiv0\pmod{2}$. This implies that
\[\sum_{i=1}^n(2\alpha_{i,2}\beta_{i,2}(2(\alpha_{i,2}^{~~~2}-\beta_{i,2}^{~~~2})+3(\beta_{i,3}-\alpha_{i,3}))+2(\alpha_{i,4}\beta_{i,2}-\alpha_{i,2}\beta_{i,4})-(\alpha_{i,3}\beta_{i,2}^{~~~2}-\beta_{i,3}\alpha_{i,2}^{~~~2})\equiv0\pmod{2}.\]
This is a contradiction. We conclude that $\pi(x+x^5)\neq0$. Therefore, $|\pi(x+x^5)|=2$.

\item Assume that $\pi(x+x^7)\in\langle\pi(x+x^5)\rangle$. This implies that there exists $\epsilon\in\{0,1\}$ such that $\pi(x+x^7)=\pi(x+x^5)^\epsilon$. 
This holds if and only if there exists $f\in [\Gamma_2,\Gamma_2]$ such that 
\[p'(x+x^7+O(x^8))=p'(x+\epsilon x^5+O(x^8))\circ f.\]
This is equivalent to 
 \begin{align*}
    p'(x+x^7+O(x^8))=&p'\Bigg(x+\left(\sum_{i=1}^n(\alpha_{i,2}\beta_{i,2}(\beta_{i,2}-\alpha_{i,2})+(\alpha_{i,3}\beta_{i,2}-\alpha_{i,2}\beta_{i,3}))\right)x^4\\
    +&\Bigg(\sum_{i=1}^n(2\alpha_{i,2}\beta_{i,2}(2(\alpha_{i,2}^{~~~2}-\beta_{i,2}^{~~~2})+3(\beta_{i,3}-\alpha_{i,3}))\\
    +&2(\alpha_{i,4}\beta_{i,2}-\alpha_{i,2}\beta_{i,4})-(\alpha_{i,3}\beta_{i,2}^{~~~2}-\beta_{i,3}\alpha_{i,2}^{~~~2})+\epsilon\Bigg)x^5\\
    +&\left(\sum_{i=1}^n(\alpha_{i,3}\beta_{i,2}-\alpha_{i,2}\beta_{i,3})\right)x^7+O(x^8)\Bigg).
\end{align*}

On the one hand, note that this implies that 
\[\sum_{i=1}^n(\alpha_{i,3}\beta_{i,2}-\alpha_{i,2}\beta_{i,3})\equiv1\pmod{2}.\]

On the other hand, note that the latter equation implies that 
\[\sum_{i=1}^n(\alpha_{i,2}\beta_{i,2}(\beta_{i,2}-\alpha_{i,2})+(\alpha_{i,3}\beta_{i,2}-\alpha_{i,2}\beta_{i,3}))=0.\]

Since $ab(b-a)\equiv0\pmod{2}$ for all $a,b\in\Z$, the latter equation implies that
\[\sum_{i=1}^n(\alpha_{i,3}\beta_{i,2}-\alpha_{i,2}\beta_{i,3})\equiv0\pmod{2}.\]
This is a contradiction. We conclude that $\pi(x+x^7)\not\in\langle\pi(x+x^5)\rangle$.
\end{enumerate}\vspace{-0.735cm}
\end{proof}

The proof of the following proposition is analogous to the general case, thus we omit it. 

\prop $T(H_1(\J_2(\Z)))\cong(\Z/2\Z)^2$.\\

By analogous methods to the general case we can prove the following statements:

\prop Let $f,g\in\Gamma_4$, with $f=p'(x+\sum_{i=4}^15\alpha_ix^i+O(x^{16}))$ and $g=p'(x+\sum_{i=4}^15\beta_ix^i+O(x^{16})$. Then
\begin{align*}
    f\circ &g=p'(x+(\alpha_4+\beta_4)x^4+(\alpha_5+\beta_5)x^5+(\alpha_6+\beta_6)x^6+(\alpha_7+\beta_7+4\alpha_4\beta_4)x^7+(\alpha_8+\beta_8+5\alpha_5\beta_4+4\alpha_4\beta_5)x^8\\
    +&(\alpha_9+\beta_9+6\alpha_6\beta_4+5\alpha_5\beta_5+4\alpha_4\beta_6)x^9+(\alpha_{10}+\beta_{10}+7\alpha_7\beta_4+6\alpha_6\beta_5+5\alpha_5\beta_6+4\alpha_4\beta_7+6\alpha_4\beta_4^{~2})x^{10}\\
    +&(\alpha_{11}+\beta_{11}+8\alpha_8\beta_4+7\alpha_7\beta_5+6\alpha_6\beta_6+5\alpha_5\beta_7+4\alpha_4\beta_8+10\alpha_5\beta_4^{~2}+12\alpha_4\beta_4\beta_5)x^{11}\\
    +&(\alpha_{13}+\beta_{13}+\alpha_9\beta_5+\alpha_7\beta_7+\alpha_5\beta_9+\alpha_7\beta_4)x^{13}+(\alpha_{15}+\beta_{15}+\alpha_{11}\beta_5+\alpha_9\beta_7+\alpha_7\beta_9+\alpha_5\beta_{11}+\alpha_7\beta_5)x^{15}+O(x^{16})).
\end{align*}

\prop Let $f\in\Gamma_4$, with $f=p'(x+\sum_{i=4}^15\alpha_ix^i+O(x^{16}))$. Then
\begin{align*}
    f^{-1}=p'(x&+(-\alpha_4)x^4+(-\alpha_5)x^5+(-\alpha_6)x^6+(-\alpha_7+4\alpha_4^{~2})x^7+(-\alpha_8+9\alpha_5\alpha_4)x^8\\
    +&(-\alpha_9+10\alpha_6\alpha_4+5\alpha_5^{~2})x^9+(-\alpha_{10}+11\alpha_7\alpha_4+11\alpha_6\alpha_5-22\alpha_4^{~3})x^{10}\\
    +&(\alpha_{11}+12\alpha_8\alpha_4+12\alpha_7\alpha_5+6\alpha_6^{~2}-78\alpha_5\beta_4^{~2})x^{11}\\
    +&(\alpha_{13}+\alpha_7+\alpha_5+\alpha_7\beta_4)x^{13}+\alpha_{15}x^{15}+O(x^{16})).
\end{align*}

\cor Let $f,g\in\Gamma_4$, with $f=p'(x+\sum_{i=4}^15\alpha_ix^i+O(x^{16}))$ and $g=p'(x+\sum_{i=4}^15\beta_ix^i+O(x^{16})$. Then
\begin{align*}
    [f,g]=p'(x&+(\alpha_5\beta_4-\beta_5\alpha_4)x^8
    +2(\alpha_6\beta_4-\beta_6\alpha_4)x^9+(3(\alpha_7\beta_4-\beta_7\alpha_4)+(\alpha_6\beta_5-\beta_6\alpha_5)+6\alpha_4\beta_4(\beta_4-\alpha_4))x^{10}\\
    +&(4(\alpha_8\beta_4-\beta_8\alpha_4)+2(\alpha_7\beta_5-\beta_7\alpha_5)+10(\alpha_5\beta_4^{~2}-\beta_5\alpha_4^{~2})+12\alpha_4\beta_4(\beta_5-\alpha_5))x^{11}\\
    +&(\alpha_7\beta_4-\beta_7\alpha_4)x^{13}+(\alpha_7\beta_5-\beta_7\alpha_5)x^{15}+O(x^{16})).
\end{align*}

From here the proof of the following proposition is analogous to the case for $k=2$.

\prop $T(H_1(\J_4(\Z)))\cong\Z/4\Z\bigoplus\Z/2\Z.$

\subsubsection{The case $k=3.$}

As in the latter subsection, by analogous methods to the general case we prove the following statements:

\prop Let $f,g\in\Gamma_3$, with $f=p'(x+\sum_{i=3}^9\alpha_ix^i+O(x^{10}))$ and $g=p'(x+\sum_{i=2}^9\beta_ix^i+O(x^{10})$. Then
\begin{align*}
    f\circ g=p'(x&+(\alpha_3+\beta_4)x^3+(\alpha_4+\beta_4)x^4+(\alpha_5+\beta_5+3\alpha_3\beta_3)x^5+(\alpha_6+\beta_6+4\alpha_4\beta_3+3\alpha_3\beta_4)x^6\\
    +&(\alpha_7+\beta_7+\alpha_5\beta_3+\alpha_3\beta_3+\alpha_3\beta_5)x^7+(\alpha_9+\beta_9+\alpha_7\beta_3+\alpha_5\beta_5+\alpha_3\beta_7+\alpha_3\beta_3+\alpha_3\beta_4)x^9+O(x^{10})).
\end{align*}

\cor Let $f\in\Gamma_3$, with $f=p'(x+\sum_{i=3}^9\alpha_ix^i+O(x^{10}))$. Then
\begin{align*}
    f^{-1}=p'(x&+(-\alpha_3)x^3+(-\alpha_4)x^4+(-\alpha_5+3\alpha_3^{~2})x^5+(-\alpha_6+4\alpha_4\alpha_3+3\alpha_3\alpha_4)x^6\\
    +&\alpha_7x^7+(\alpha_9+\alpha_3+\alpha_5+\alpha_3\alpha_4+\alpha_3\alpha_5)x^9+O(x^{10})).
\end{align*}

\cor Let $f,g\in\Gamma_3$, with $f=p'(x+\sum_{i=3}^9\alpha_ix^i+O(x^{10}))$ and $g=p'(x+\sum_{i=2}^9\beta_ix^i+O(x^{10})$. Then
\[[f,g]=p'((\alpha_4\beta_3-\beta_4\alpha_3)x^6+(\alpha_4\beta_3-\beta_4\alpha_3)x^9+O(x^{10}).\]

The proof of the following statement is analogous to the general case:

\prop The isomorphism $T(H_1(\J_3(\Z)))\cong(\Z/2\Z)^2$ holds.

\subsection{The free submodule of $H_1(\J_k(R))$}

In order to complete our computation of $H_1(\J_k(\Z))$, it remains to determine the free submodule of this group. We can go back to generalize some of the propositions of this section.
\\

Let $R$ be a domain. We know that $[\J_k(R),\J_k(R)]\subset \J_{2k}(R)$. By the factor theorem, this implies that the following diagram commutes:

 \begin{center}
        \begin{tikzcd}
            \J_k(R)\arrow{rr}{p}\arrow{ddr}{\pi}&&\J_k^{2k}(R)\\
            &&&\\
            &H_1(\J_k(R))\arrow{uur}{\varphi}&
        \end{tikzcd}
\end{center}
Since $p$ is surjective, the morphism $\varphi$ is surjective. From the diagram we can deduce that $ker(\varphi)=\pi(\J_{2k}(R))$. From  here, we have two  ways forward.

\prop If $R$ is a domain such that $2\in R^*$, then 
\[H_1(\J_k(R))\cong \J_k^{2k}(R).\]
\begin{proof}
    By corollary \ref{cor:dosss}, we have that $[\J_k(R),\J_k(R)]=\J_{2k}(R)$. This implies that $\pi(\J_{2k}(R))=\{0\}$. We conclude that $\varphi$ is an isomorphism.
\end{proof}

\prop If $R$ is a torsion-free domain, then 
 \[0\to T(H_1(\J_k(R)))\to H_1(\J_k(R))\overset{\varphi}{\to}\J_k^{2k}(R)\to 0\]
 is a group extension.
 \\
 
 For our case of interest, it suffices to note that 
 \[\{\pi(x+x^i)|k\leqslant i\leqslant  2k-1\}\cup\{\pi(x+x^{2i+1})|k\leqslant i\leqslant \frac{c_k}{2}-1\}\]
 is a finite set of generators for $H_1(\J_k(\Z))$. This implies that $H_1(\J_k(\Z))$ is a finitely generated abelian group, which in turn implies that 
 \[H_1(\J_k(\Z))\cong \J_k^{2k}(\Z)\bigoplus \pi(\J_k(\Z)).\]
 For this reason, we study the structure of $\J_k^{2k}(R)$.

\prop The group $\J_k^{2k}(R)$ is isomorphic to a subgroup of $R^k$.
\begin{proof}
Define 
 \begin{align*}
     \psi:\J_k^{2k}(R)&\to R^k\\
     x+\sum_{i=k}^{2k-1}\alpha_ix^i+O(x^{2k})&\mapsto(\alpha_i)_{i=k}^{2k-2}\oplus( k\alpha_k^{~2}-2\alpha_{2k-1}).
 \end{align*}
In order to prove that this is a group morphism, we will describe the product in $\J_k^{2k}(R)$. Let $f$ and $g$ be elements in $\J_k^{2k}(R)$, with $f=x+\sum_{i=k}^{2k-1}\alpha_ix^i+O(x^{2k-1})$ and $g=x+\sum_{i=k}^{2k-1}\beta_ix^i+O(x^{2k})$. By Proposition \ref{prop:formulas} we have that
\[f\circ g=x+\sum_{i=k}^{2k-2}(\alpha_i+\beta_i)x^i+(\alpha_{2k-1}+\beta_{2k-1}+k\alpha_k\beta_k)x^{2k-1}+O(x^{2k}).\]

Since all but the last coordinate of $\psi$ is a projection, it suffices to note that 
    \[k(\alpha_k+\beta_k)^2-2(\alpha_{2k-1}+\beta_{2k-1}+k\alpha_k\beta_k)=k\alpha_k^{~2}-2\alpha_{2k-1}+k\beta_k^{~2}-2\beta_{2k-1}.\]
We conclude that $\psi$ is a group morphism.\\

In order to prove that $\psi$ is injective, note that $\psi(x+\sum_{i=k}^{2k-1}\alpha_ix^i)=0$ if and only if $2\alpha_{2k-1}=k\alpha_k^{~2}$ and $\alpha_i=0$ for all $i\in\{k,\dots,2k-2\}$. The second condition implies that $\alpha_k=0$, thus $\alpha_{2k-1}=0$. We conclude that $ker(\psi)=\{x+O(x^{2k})\}$. It follows that $\J_k^{2k}(R)\cong Im(\psi)<R^n$.
\end{proof}

We take this paragraph to note that the last coordinate in the function $\psi$ is known in the literature as the residue of order $k$ of $x+\sum_{i\geqslant k}\alpha_ix^i$. More about them and their applications can be found in \cite{Eynard-Bontemps_Navas_2023}.
\\

We conclude with the following corollary:

\cor The isomorphism $\J_k^{2k}(\Z)\cong \Z^k$ holds.

\begin{proof}
We will prove that $[\Z^k:Im(\psi)]<\infty$, since this will imply the isomorphism. Let $(n_1,\dots,n_k)\in\Z^k$. Note that $n_k=kn_1^{~2}-\Tilde{n}$. Moreover, there exists $n\in\Z$ and $\epsilon\in\{0,-1\}$ such that $\Tilde{n}=2n+\epsilon$. Thus $(n_1,\dots n_{k-1},kn_1^{~2}-2n)+(0,\dots,0,-\epsilon)$. This implies that 
\[\Z^k/\psi(\J_k^{2k}(\Z))\cong\Z/2\Z.\]

We conclude that $[\Z^k:Im(\psi)]=2<\infty$.

\end{proof}

Finally, it suffices to join the latter proposition with Theorem \ref{thm:tors} in order to prove Theorem \ref{thm:princ}.

\section{Appendix: The product on a quotient}
Let $R$ be a torsion-free domain. In this section we will give a description for the product in $\J_k^{c_k}(R)$, for $k\geqslant5$. Later, we give a description for the product, inverse and commutator in $\Gamma_k$. These descriptions can be found at the end 
of this section.
\\

Let $f$ and $g$ be elements in $\J_k^{c_k}(R)$, with $f=x+\sum_{i=k}^{c_k-1}\alpha_ix^i+O(x^{c_k})$ and $g=x+\sum_{i=k}^{c_k-1}\beta_ix^i+O(x^{c_k})$. We know that $f\circ g=x+\sum_{i=k}^{c_k-1}\beta_ix^i+\sum_{i=k}^{c_k-1}\alpha_i(x+\sum_{j=k}^{c_k-1}\beta_jx^j)^i+O(x^{c_k})$. We would like to write the coefficients of $f\circ g$ in function of the coefficients of $f$ and $g$. For this purpose, we will make use of the multinomial theorem. In order to state this theorem comfortably, we introduce some definitions.

\deff Let $v\in (\Z_{\geqslant0})^m$ and $i\in\{1,\dots,m\}$. We write $v_i$ for the $i^{th}$ coordinate of $v$.

\deff For $v\in (\Z_{\geqslant0})^m$ we define $|v|:=\sum_{i=1}^mv_i$.

\deff Given $v\in(\Z_{\geqslant0})^m$ such that $|v|=n$, we define the number
\[\binom{n}{v}:=\frac{n!}{\prod_{i=1}^nv_i!}.\]

\begin{thm}[Multinomial Theorem.]
    For a finite collection $\{\alpha_i\}_{i=1}^m$ in $R$ one has that
    \[\left(\sum_{i=1}^m\alpha_i\right)^n=\sum_{\substack{v\in(\Z_{\geqslant0})^m\\ |v|=n}}\binom{n}{v}\prod_{i=1}^m(\alpha_i)^{v_i}.\]
\end{thm}

In order to fit this theorem to our case we will introduce some notation. Note that an element of $\J_k^{c_k}(R)$ has $c_k-k+1$ terms. We define $m_k:=c_k-k+1$. We fix $w^k\in(\Z_{\geqslant0})^{m_k}$ as the tuple of coordinates $(w^k)_1=1$ and $(w^k)_i=k+(i-2)$ for $i\in\{2,\dots,m_k\}$. For $v$ and $w$ in $(\Z_{\geqslant0})^m$, we define $v\cdot w:=\sum_{i=1}^mv_iw_i$. Since any element of $\J_k^{c_k}(R)$ is in particular an element of $R[[x]]/x^{c_k}R[[x]]$, any term of $f\circ g$ that is divisible by $x^c_k$ vanishes. Thus, we define 
\[S_{n,k}:=\{v\in(\Z_{\geqslant0})^{m_k}|v\cdot w^k<c_k\text{ and }|v|=n\}.\]

We can write $f\circ g$ as follows:

\cor \label{cor:for} Let $f$ and $g$ be elements in $\J_k^{c_k}(R)$, with $f=x+\sum_{i=k}^{c_k-1}\alpha_ix^i+O(x^{c_k})$ and $g=x+\sum_{i=k}^{c_k-1}\beta_ix^i+O(x^{c_k})$. Then
\[f\circ g=x+\sum_{i=k}^{c_k-1}\beta_ix^i+\sum_{l=k}^{c_k-1}\alpha_l\sum_{v\in S_{l,k}}\binom{l}{v}\left(\prod_{i=k}^{c_k-1}(\beta_i)^{v_{i-k+2}}\right)x^{v\cdot w^k}+O(x^{c_k}).\]

\begin{proof}
    Note that 
    \[f\circ g=x+\sum_{i=k}^{c_k-1}\beta_ix^i+\sum_{i=k}^{c_k-1}\alpha_i\left(x+\sum_{j=k}^{c_k-1}\beta_jx^j\right)^i+O(x^{c_k}).\]
  Let $l\geqslant k$. By the multinomial theorem we have that
    \begin{align*}
        \left(x+\sum_{j=k}^{c_k-1}\beta_jx^j\right)^l=&\sum_{\substack{v\in(\Z_{0\geqslant})^{m_k}\\ |v|=l}}\binom{l}{v}x^{v_1}\prod_{j=k}^{c_k-1}(\beta_jx^j)^{v_{j-k+2}}\\
        % =&\sum_{\substack{v\in(\Z_{0\geqslant})^{m_k}\\ |v|=l}}\binom{l}{v}\prod_{j=k}^{c_k-1}(\beta_j)^{v_{j-k+2}}x^{v_1}\prod_{j=k}^{c_k-1}(x^j)^{v_{j-k+2}}\\
        =&\sum_{\substack{v\in(\Z_{0\geqslant})^{m_k}\\ |v|=l}}\binom{l}{v}\prod_{j=k}^{c_k-1}(\beta_j)^{v_{j-k+2}}x^{v_1}\prod_{i=2}^{m_k}(x^{k+(i-2)})^{v_i}\\
        % =&\sum_{\substack{v\in(\Z_{0\geqslant})^{m_k}\\ |v|=l}}\binom{l}{v}\prod_{j=k}^{c_k-1}(\beta_j)^{v_{j-k+2}}x^{v_1}\prod_{i=2}^{m_k}(x^{(w^k)_i})^{v_i}\\
        % =&\sum_{\substack{v\in(\Z_{0\geqslant})^{m_k}\\ |v|=l}}\binom{l}{v}\prod_{j=k}^{c_k-1}(\beta_j)^{v_{j-k+2}}x^{v_1}\prod_{i=2}^{m_k}x^{(w^k)_iv_i}\\
        =&\sum_{\substack{v\in(\Z_{0\geqslant})^{m_k}\\ |v|=l}}\binom{l}{v}\prod_{j=k}^{c_k-1}(\beta_j)^{v_{j-k+2}}\prod_{i=1}^{m_k}x^{(w^k)_iv_i}\\
        % =&\sum_{\substack{v\in(\Z_{0\geqslant})^{m_k}\\ |v|=l}}\binom{l}{v}\prod_{j=k}^{c_k-1}(\beta_j)^{v_{j-k+2}}x^{\sum_{i=1}^{m_k}(w^k)_iv_i}
        =&\sum_{\substack{v\in(\Z_{0\geqslant})^{m_k}\\ |v|=l}}\binom{l}{v}\prod_{j=k}^{c_k-1}(\beta_j)^{v_{j-k+2}}x^{w^k\cdot v}
    \end{align*}
    If $v\in(\Z_{\geqslant0})^{m_k}$ is such that $|v|=l$ and $v\cdot w^k\geqslant c_k$, then
    \[\binom{l}{v}\prod_{j=k}^{c_k-1}(\beta_j)^{v_{j-k+2}}x^{w^k\cdot v}\equiv0\pmod{x^{c_k}}.\]
    On the other hand, if $v\in(\Z_{\geqslant0})^{m_k}$ is such that $|v|=l$ and
    \[\binom{l}{v}\prod_{j=k}^{c_k-1}(\beta_j)^{v_{j-k+2}}x^{w^k\cdot v}\not\equiv0\pmod{x^{c_k}},\]
    then $v\cdot w^k<c_k$. Thus, $v\in S_{l,k}$. We conclude that 
    \[\sum_{\substack{v\in(\Z_{\geqslant0})^{m_k}\\ |v|=l}}\binom{l}{v}\prod_{j=k}^{c_k-1}(\beta_j)^{v_{j-k+2}}x^{w^k\cdot v}=\sum_{v\in S_{l,k}}\binom{l}{v}\prod_{j=k}^{c_k-1}(\beta_j)^{v_{j-k+2}}x^{w^k\cdot v}.\]
\end{proof}

We will study the sets $S_{l,k}$, in order to figure out which monomials appear as terms of $f\circ g$.

\lem Let $k$ and $l$ be integers such that $l\geqslant k\geqslant2$. If $v\in S_{l,k}$, then 
\[\#\{i\in\{2,\dots,m_k\}|v_i\neq0\}\leqslant2.\]
\begin{proof}
 Let $v\in S_{l,k}$. Assume that there exists $i_1,i_2,i_3\in\{2,\dots,m_k\}$ such that $v_{i_1},v_{i_2},v_{i_3}$ are non zero. Without loss of generality, assume that $i_1<i_2<i_3$. If $\max\{v_1,v_2,v_3\}\geqslant2$, then
 \begin{align*}
     v\cdot w^k&\geqslant v_{i_1}(k+(i_1-2))+v_{i_2}(k+(i_2-2))+v_{i_3}(k+(i_3-2))\\ &\geqslant v_{i_1}k+v_{i_2}(k+1)+v_{i_3}(k+2)\\
&     \geqslant  2k+k+1+k+2=4k+3\\
&     \geqslant c_k.
 \end{align*}
 This contradicts that $v\in S_{l,k}$. Thus $v_{i_1}=v_{i_2}=v_{i_3}=1$. \\
 
 Assume that there exists $i_4\in\{2,\dots,m_k\}\setminus \{i_1,i_2,i_3\}$ such that $v_{i_4}\neq0$. Analogously to the previous case, we conclude that $v_{i_4}=1$. Without loss of generality, assume that $i_3<i_4$. Note that
  \begin{align*}
      v\cdot w^k&\geqslant v_{i_1}(k+(i_1-2))+v_{i_2}(k+(i_2-2))+v_{i_3}(k+(i_3-2))+v_{i_4}(k+(i_4-2)) \\
     & =k+(i_1-2)+k+(i_2-2)+k+(i_3-2)+k+(i_4-2) \\ 
&      \geqslant k+k+1+k+2+k+3\\
     & =4k+6\\
&      \geqslant c_k.
  \end{align*}
  This is a contradiction. We conclude that $\{i\in\{2,\dots,m_k\}|v_i\neq0\}=\{i_1,i_2,i_3\}$.\\
  
  Since $|v|=l$, we have that $v_1=l-3$. This implies that
  \begin{align*}
      v\cdot w^k&=v_1+ v_{i_1}(k+i_1-2)+v_{i_2}(k+i_2-2)+v_{i_3}(k+i_3-2) \\
  &=l-3+k+i_1-2+k+i_2-2+k+i_3-2 \\
&  \geqslant k-3+k+k+1+k+2=4k\\
&  \geqslant c_k.
  \end{align*}
 This is a contradiction. We conclude that $\#\{i\in\{2,\dots,m_k\}|v_i\neq0\}\leqslant2$.
\end{proof}

We will restrict ourselves to the case of $k\geqslant5$, since this will allow us to describe $S_{l,k}$ in a general fashion. For the three cases not included, it suffices to compute the sets $S_{l,k}$ by hand.

\lem Let $k$ and $l$ be integers such that $l\geqslant k\geqslant5$. Given $v\in S_{l,k}$, we have the following possibilities:
\begin{enumerate}[$a)$]
    \item If $\{i\in\{2,\dots,m_k\}|v_i\neq0\}=\{i_1,i_2\}$, then $v_{i_1}=1=v_{i_2}$.

    \item If $\{i\in\{2,\dots,m_k\}|v_i\neq0\}=\{i\}$, then $v_i\leqslant2$.
\end{enumerate}
\begin{proof}
    ~\begin{enumerate}[$a)$]
        \item Without loss of generality, assume that $i_1<i_2$. Since $2\leqslant i_1$, we have that $k+(i_1-2)\geqslant k$ and $k+(i_2-2)\geqslant k+1$. Assume that $v_{i_1}\geqslant2$. This implies that
        \begin{align*}
            v\cdot w^k&=v_1+v_{i_1}(k+(i_1-2))+v_{i_2}(k+(i_2-2))\\
            &=l-v_{i_1}-v_{i_2}+v_{i_1}(k+(i_1-2))+2(k+(i_2-2))\\ 
            &=l+v_{i_1}(k+(i_1-2)-1)+v_{i_2}(k+(i_2-2)-1)\\ 
            &\geqslant k+2(k-1)+k\\
            &=4k-2\\
            &\geqslant c_k.
        \end{align*}
        This is a contradiction, thus $v_{i_1}=1$. This implies that $v_1=l-1-v_{i_2}$. If we  assume that $v_{i_2}\geqslant2$, then 
        \begin{align*}
            v\cdot w^k&=v_1+v_{i_1}(k+(i_1-2))+v_{i_2}(k+(i_2-2))\\
            &=l-1-v_{i_2}+k+(i_1-2)+v_{i_2}(k+(i_2-2))\\ 
            &=l-1+k+(i_1-2)+v_{i_2}(k+(i_2-2)-1)\\ 
            &\geqslant k-1+k+2k\\
            &=4k-1\\
            &\geqslant c_k.
        \end{align*}
       This is a contradiction. We conclude that $v_{i_1}=1=v_{i_2}$.

        \item Assume that $v_i\geqslant4$. This implies that
        \[v\cdot w^k\geqslant 4k\geqslant c_k.\]
        This is a contradiction. We conclude that $v_i\leqslant3$.\\
        
        Assume that $v_i=3$. Note that $|v|=l$ implies that $v_1=l-3$. We have that
        \begin{align*}
            v\cdot w^k&=v_1+v_i(k+(i-2))\\
            &=l-3+3(k+(i-2))\\ 
            &\geqslant k-3+3k\\
            &= 4k-3\\
            &\geqslant c_k.
        \end{align*}
        This is a contradiction. We conclude that $v_i\leqslant2$.
    \end{enumerate}
\end{proof}
\obs Let $l\geqslant k$. If $v\in S_{l,k}$, then $v_1\neq0$.
\deff Let $l\geqslant k$. Let $v\in S_{l,k}$. We say that $v$ is of:
\begin{enumerate}[$a)$]
    \item Type A if it has two nonzero coordinates and $v_1=l-1$.

    \item Type B if it has two nonzero coordinates and $v_1=l-2$.

    \item Type C if it has three nonzero coordinates, in which case $v_1=l-2$.
\end{enumerate}
 
This lead us us to describe $S_{l,k}$ as the union of $\{(l,0,\dots,0)\}$ and the vectors of types A, B and C contained in it. With this, we can classify the monomials of $f\circ g$.

\lem Let $k\geqslant5$. Let $f$ and $g$ be elements in $\J_k^{c_k}(R)$, with $f=x+\sum_{i=k}^{c_k-1}\alpha_ix^i+O(x^{c_k})$ and $g=x+\sum_{i=k}^{c_k-1}\beta_ix^i+O(x^{c_k})$. Let $l\geqslant k$ and $v\in S_{l,k}$. We have that:
\begin{enumerate}[$a)$]
    \item If $v$ is of type A, with $i\geqslant k$ such that $v_{i-k+2}=1$, then
    \[\alpha_l\binom{l}{v}\left(\prod_{e=k}^{c_k-1}(\beta_e)^{v_{e-k+2}}\right)x^{v\cdot w^k}=l\alpha_l\beta_ix^{l-1+i},\]
    and $v\cdot w^k\geqslant2k-1$.

    \item If $v$ is of type B, with $i\geqslant k$ such that $v_{i-k+2}=2$, then
    \[\alpha_l\binom{l}{v}\left(\prod_{e=k}^{c_k-1}(\beta_e)^{v_{e-k+2}}\right)x^{v\cdot w^k}=\binom{l}{2}\alpha_l\beta_i^{~2}x^{l-2+2i},\]
    and $v\cdot w^k\geqslant3k-2$.

    \item If $v$ is of type C, with $j>i\geqslant k$ such that $v_{i-k+2}=1=v_{j-k+2}$, then
    \[\alpha_l\binom{l}{v}\left(\prod_{e=k}^{c_k-1}(\beta_e)^{v_{e-k+2}}\right)x^{v\cdot w^k}=l(l-1)\alpha_l\beta_i\beta_jx^{l-2+i+j},\]
    and $v\cdot w^k\geqslant3k-1$.

    \item If $v=(l,0,\dots,0)$, then
    \[\alpha_l\binom{l}{v}\left(\prod_{e=k}^{c_k-1}(\beta_e)^{v_{e-k+2}}\right)x^{v\cdot w^k}=\alpha_lx^l.\]
\end{enumerate}
\begin{proof}~
    \begin{enumerate}[$a)$]
        \item Note that $v=(l-1,0\dots,0,1,0\dots,0)$. This implies that 
        \[\binom{l}{v}=\frac{l!}{(l-1)!0!\cdots1!\cdots0!}=l\]
        and
        \[\prod_{e=k}^{c_k-1}(\beta_e)^{v_{e-k+2}}=\beta_i.\]
        
        Furthermore, we have that 
        \[v\cdot w^k=l-1+i\geqslant k-1+i\geqslant k-1+k=2k-1.\]
        
        \item Note that $v=(l-2,0\dots,0,2,0\dots,0)$. This implies that
        \[\binom{l}{v}=\frac{l!}{(l-2)!0!\cdots2!\cdots0!}=\binom{l}{2}\]
        and
        \[\prod_{e=k}^{c_k-1}(\beta_e)^{v_{e-k+2}}=\beta_i^{~2}.\]
        
        Furthermore, we have that 
        \[v\cdot w^k=l-2+2i\geqslant k-2+2i\geqslant k-2+2k=3k-2.\]

        \item Note that $v=(l-2,0\dots,0,1,0\dots,0,1,0,\dots,0)$. This implies that
        \[\binom{l}{v}=\frac{l!}{(l-2)!0!\cdots1!\cdots0!1!\dots0!}=l(l-1)\]
        and
        \[\prod_{e=k}^{c_k-1}(\beta_e)^{v_{e-k+2}}=\beta_i\beta_j.\]
        
        Furthermore, since $j>i\geqslant k$, we have that $j\geqslant k+1$. This implies that
        \[v\cdot w^k=l-2+i+j\geqslant k-2+i+j\geqslant k-2+k+k+1=3k-1.\]

        \item We have that
        \[\binom{l}{v}=\frac{l!}{l!0!\cdots0!\cdots0!}=1\]
        and 
        \[\prod_{e=k}^{c_k-1}(\beta_e)^{v_{e+k-2}}=1.\]
        
        Furthermore, we have that $v\cdot w^k=l$. This proves the equality. 
    \end{enumerate}
\end{proof}

    From here, it suffices to order the monomials of $f\circ g$ and factor the common powers of $x$ in order to prove the following.
    
\begin{thm}
    Let $k\geqslant5$. Let $f$ and $g$ be elements in $\J_k^{c_k}(R)$, with $f=x+\sum_{i=k}^{c_k-1}\alpha_ix^i+O(x^{c_k})$ and $g=x+\sum_{i=k}^{c_k-1}\beta_ix^i+O(x^{c_k})$. Let $\gamma_l$ be the $l$th coefficient of $f\circ g$.
    \begin{enumerate}[$a)$]
        \item If $k\leqslant l\leqslant2k-2$, then $\gamma_l=\alpha_l+\beta_l$.

        \item If $2k-1\leqslant l\leqslant 3k-3$, then $\gamma_l=\alpha_l+\beta_l+\sum_{i=k}^{l-k+1}i\alpha_i\beta_{l-i+1}$.

        \item $\gamma_{3k-2}=\alpha_{3k-2}+\beta_{3k-2}+\sum_{i=k}^{2k-1}i\alpha_i\beta_{3k-1-i}$+$\binom{k}{2}\alpha_k\beta_k^{~2}$.

        \item $\gamma_{3k-1}=\alpha_{3k-1}+\beta_{3k-1}+\sum_{i=k}^{2k}i\alpha_i\beta_{3k-i}$+$\binom{k+1}{2}\alpha_{k+1}\beta_k^{~2}+k(k-1)\alpha_k\beta_k\beta_{k+1}$.

        \item \[\gamma_{3k}=\alpha_{3k}+\beta_{3k}+\sum_{i=k}^{2k+1}i\alpha_i\beta_{3k+1-i}+\binom{k}{2}\alpha_k\beta_{k+1}^{~2}+\binom{k+2}{2}\alpha_{k+2}\beta_k^{~2}+k(k-1)\alpha_k\beta_k\beta_{k+2}+(k+1)k\alpha_{k+1}\beta_k\beta_{k+1}.\]

        \item If $k\equiv0\pmod{2}$, then
        \[\gamma_{3k+1}=\alpha_{3k+1}+\beta_{3k+1}+\sum_{i=k}^{2k+2}i\alpha_i\beta_{3k+2-i}+\binom{k+1}{2}\alpha_{k+1}(\beta_{k+1})^2+\binom{k+3}{2}\alpha_{k+3}\beta_k^{~2}+\sum_{\substack{k\leqslant l,i,j<c_k\\
     i<j\\l-2+i+j=3k+1}}l(l-1)\alpha_l\beta_i\beta_j.\]

     \item If $k\equiv0\pmod{4}$, then
        \begin{align*}
            \gamma_{3k+2}=&\alpha_{3k+2}+\beta_{3k+2}+\sum_{i=k}^{2k+3}i\alpha_i\beta_{3k+3-i}+\binom{k}{2}\alpha_k(\beta_{k+2})^2+\binom{k+2}{2}\alpha_{k+2}(\beta_{k+1})^2+\binom{k+4}{2}\alpha_{k+4}\beta_k^{~2}\\
            +&\sum_{\substack{k\leqslant l,i,j<c_k\\
     i<j\\l-2+i+j=3k+2}}l(l-1)\alpha_l\beta_i\beta_j.
        \end{align*}
and
        \begin{align*}          
        \gamma_{3k+3}=&\alpha_{3k+3}+\beta_{3k+3}+\sum_{i=k}^{2k+4}i\alpha_i\beta_{3k+4-i}+\binom{k+1}{2}\alpha_{k+1}(\beta_{k+2})^2+\binom{k+3}{2}\alpha_{k+3}(\beta_{k+1})^2+\binom{k+5}{2}\alpha_{k+5}\beta_k^{~2}\\
            +&\sum_{\substack{k\leqslant l,i,j<c_k\\
     i<j\\l-2+i+j=3k+3}}l(l-1)\alpha_l\beta_i\beta_j.\end{align*}
    \end{enumerate}
    \end{thm}

\cor Let $k\geqslant5$ Let $f$ and $g$ be elements in$\Gamma_k$, with $f=p'(x+\sum_{i=k}^{c_k-1}\alpha_ix^i+O(x^{c_k}))$ and $g=p'(x+\sum_{i=k}^{c_k-1}\beta_ix^i+O(x^{c_k}))$. Let $\gamma_l$ be the $l$th coefficient of $f\circ g$.
    \begin{enumerate}[$a)$]
        \item If $k\leqslant l\leqslant2k-2$, then $\gamma_l=\alpha_l+\beta_l$.

        \item If $2k-1\leqslant l\leqslant d_k-1$, then $\gamma_l=\alpha_l+\beta_l+\sum_{i=k}^{l-k+1}i\alpha_i\beta_{l-i+1}$.

        \item If $d_k\leqslant l\leqslant 3k-3$, then 
        
        \[\gamma_l=\left\{\begin{array}{cc}
            \alpha_l+\beta_l+\sum_{i=k}^{l-k+1}i\alpha_i\beta_{l-i+1}\pmod{2} &, \text{ if } l\equiv1\pmod{2}  \\
             0 &, \text{ if } l\equiv0\pmod{2}
        \end{array}\right.\].

        \item If $k\equiv0\pmod{2}$, then $\gamma_{3k-2}=0$. If $k\equiv1\pmod{2}$, then 
        \[\gamma_{3k-2}\equiv\alpha_{3k-2}+\beta_{3k-2}+\sum_{i=k}^{2k-1}i\alpha_i\beta_{3k-1-i}+\binom{k}{2}\alpha_k\beta_k^{~2}\pmod{2}.\]

        \item If $k\equiv0\pmod{2}$, then
        \[\gamma_{3k-1}\equiv\alpha_{3k-1}+\beta_{3k-1}+\sum_{i=k}^{2k}i\alpha_i\beta_{3k-i}+\binom{k+1}{2}\alpha_{k+1}\beta_k^{~2}\pmod{2}.\]
        
        If $k\equiv1\pmod{2}$, then $\gamma_{3k-1}=0$.

        \item If $k\equiv0\pmod{2}$, then $\gamma_{3k}=0$. If $k\equiv1\pmod{2}$, then
        \[\gamma_{3k}\equiv\alpha_{3k}+\beta_{3k}+\sum_{i=k}^{2k+1}i\alpha_i\beta_{3k+1-i}+\binom{k}{2}\alpha_k(\beta_{k+1})^2+\binom{k+2}{2}\alpha_{k+2}\beta_k^{~2}\pmod{2}.\]

        \item If $k\equiv0\pmod{2}$, then
        \[\gamma_{3k+1}\equiv\alpha_{3k+1}+\beta_{3k+1}+\sum_{i=k}^{2k+2}i\alpha_i\beta_{3k+2-i}+\binom{k+1}{2}\alpha_{k+1}(\beta_{k+1})^2+\binom{k+3}{2}\alpha_{k+3}\beta_k^{~2}\pmod{2}.\]

     \item If $k\equiv0\pmod{4}$, then $\gamma_{3k+2}=0$, and    
        \[\gamma_{3k+3}\equiv\alpha_{3k+3}+\beta_{3k+3}+\sum_{i=k}^{2k+4}i\alpha_i\beta_{3k+4-i}+\alpha_{k+3}(\beta_{k+1})^2\pmod{2}.\]
    \end{enumerate}
    
    With the previous formula, we can compute inverses and commutators in $\Gamma_k$. We omit the following proof, since this essentially consists of evaluations and direct computations.
    
\cor Let $k\geqslant5$. Let $f\in\Gamma_k$, with $f=p'(x+\sum_{i=k}^{c_k-1}\alpha_ix^i+O(x^{c_k}))$. Let $\beta_l$ be the $l$th coefficient of $f^{-1}$.
    \begin{enumerate}[$a)$]
        \item If $k\leqslant l\leqslant2k-2$, then $\beta_l=-\alpha_l$.

        \item If $2k-1\leqslant l\leqslant d_k-1$, then $\beta_l=-\alpha_l+\sum_{i=k}^{l-k+1}i\alpha_i\alpha_{l-i+1}$.

         \item If $d_k-1\leqslant l\leqslant 3k-3$, then

         \[\beta_l=\left\{\begin{array}{cc}
             -\alpha_l+\sum_{i=k}^{l-k+1}i\alpha_i\alpha_{l-i+1} \pmod{2} &,\text{ if } l\equiv1\pmod{2} \\
             0 &,\text{ if } l\equiv0\pmod{2} 
         \end{array}\right..\]

        \item If $k\equiv0\pmod{2}$, then $\beta_{3k-2}=0$. If $k\equiv1\pmod{2}$, then 
        \[\beta_{3k-2}\equiv\alpha_{3k-2}+\sum_{i=k}^{2k-1}i\alpha_i\alpha_{3k-1-i}+\left(\binom{k}{2}+1\right)\alpha_k^{~3}\pmod{2}.\]

        \item If $k\equiv0\pmod{2}$, then
        \[\beta_{3k-1}\equiv\alpha_{3k-1}+\sum_{i=k}^{2k}i\alpha_i\alpha_{3k-i}+\binom{k+1}{2}\alpha_{k+1}\alpha_k^{~2}\pmod{2}.\]
        If $k\equiv1\pmod{2}$, then $\beta_{3k-1}=0$.

        \item If $k\equiv0\pmod{2}$, then $\beta_{3k}=0$. If $k\equiv1\pmod{2}$, then
        \[\beta_{3k}\equiv\alpha_{3k}+\sum_{i=k}^{2k+1}i\alpha_i\alpha_{3k+1-i}+\binom{k}{2}\alpha_k(\alpha_{k+1})^2+\left(\binom{k+2}{2}+1\right)\alpha_{k+2}\alpha_k^{~2}\pmod{2}.\]

        \item If $k\equiv0\pmod{2}$, then
        \[\beta_{3k+1}\equiv\alpha_{3k+1}+\sum_{i=k}^{2k+2}i\alpha_i\alpha_{3k+2-i}+\left(\binom{k+1}{2}+1\right)(\alpha_{k+1})^3+\binom{k+3}{2}\alpha_{k+3}\alpha_k^{~2}\pmod{2}.\]

     \item If $k\equiv0\pmod{4}$, then $\beta_{3k+2}=0$, and
        \[\beta_{3k+3}\equiv\alpha_{3k+3}+\sum_{i=k}^{2k+4}i\alpha_i\alpha_{3k+4-i}\pmod{2}.\]
    \end{enumerate}
    
\cor \label{cor:conmu} Let $k\geqslant5$. Let $f$ and $g$ be elements in$\Gamma_k$, with $f=p'(x+\sum_{i=k}^{c_k-1}\alpha_ix^i+O(x^{c_k}))$ and $g=p'(x+\sum_{i=k}^{c_k-1}\beta_ix^i+O(x^{c_k}))$. Let $\delta_l$ be the $l$th coefficient of $[f,g]$.
    \begin{enumerate}
        \item If $k\leqslant l\leqslant2k-1$, then $\delta_l=0$.

        \item If $2k \leqslant l\leqslant d_k-1$, then
        \[\delta_l=\sum_{i=k}^{l-k+1}(2i-l-1)\alpha_i\beta_{l-i+1}.\]

        \item If $d_k\leqslant l\leqslant 3k-3$, then $\delta_l\equiv0\pmod{2}$. 

        \item If $k\equiv0\pmod{2}$, then $\delta_{3k-2}=0$. If $k\equiv1\pmod{2}$, then 
        \[\delta_{3k-2}\equiv\binom{k}{2}(\alpha_k\beta_k^{~2}+\beta_k\alpha_k^{~2}) \pmod{2}.\]

        \item If $k\equiv0\pmod{2}$, then \[\delta_{3k-1}\equiv\binom{k+1}{2}(\alpha_{k+1}\beta_k^{~2}+\beta_{k+1}\alpha_k^{~2})\pmod{2}.\]
        If $k\equiv1\pmod{2}$, then $\delta_{3k-1}=0$.

        \item If $k\equiv0\pmod{2}$, then $\delta_{3k}=0$. If $k\equiv1\pmod{2}$, then
        \[\delta_{3k}\equiv\binom{k}{2}(\alpha_k(\beta_{k+1})^2+\beta_k(\alpha_{k+1})^2)+\binom{k+2}{2}(\alpha_{k+2}\beta_k^{~2}+\beta_{k+2}\alpha_k^{~2})\pmod{2}.\]
        
        \item If $k\equiv0\pmod{2}$, then
        \[\delta_{3k+1}\equiv\binom{k+3}{2}(\alpha_{k+3}\beta_k^{~2}+\beta_{k+3}\alpha_k^{~2})\pmod{2}.\]

     \item If $k\equiv0\pmod{4}$, then $\delta_{3k+2}=0$.

     \item If $k\equiv0\pmod{4}$, then
        \[\delta_{3k+3}\equiv\alpha_{k+3}(\beta_{k+1})^2+\beta_{k+3}(\alpha_{k+1})^2\pmod{2}.\]
    \end{enumerate}
\bibliography{refart}
\bibliographystyle{alpha}
\nocite{*}

{\bf Acknowledgments.} I'm strongly indebted to Luis Arenas for the proof 
of Proposition 2.1.1.\\ This work has been partially
supported by FONDECYT 1210155.

\end{document}